\documentclass[11pt]{article}
\usepackage[a4paper,margin=26mm]{geometry}
\usepackage{amsmath,amssymb,amsthm,mathtools}
\usepackage[T1]{fontenc}
\usepackage{lmodern,microtype}
\usepackage[colorlinks=true,linkcolor=blue,urlcolor=blue,citecolor=blue]{hyperref}
\newtheorem{theorem}{Theorem}[section]
\newtheorem{lemma}[theorem]{Lemma}
\newtheorem{proposition}[theorem]{Proposition}
\newtheorem{corollary}[theorem]{Corollary}
\theoremstyle{remark}
\newtheorem{remark}[theorem]{Remark}
\newcommand{\E}{\mathbb E}
\newcommand{\PP}{\mathbb P}
\newcommand{\inj}{\mathrm{inj}}
\newcommand{\norm}[1]{\left\lVert#1\right\rVert}
\newcommand{\one}{\mathbf1}
\DeclareMathOperator{\Var}{Var}
\numberwithin{equation}{section}
\title{Optimal fiber comparisons for random tensors}
\author{Ziqian Lin\thanks{Corresponding author. School of Mathematics and
Statistics, Fuzhou University, Fuzhou, China.
Email \href{mailto:linziqian@fzu.edu.cn}{\texttt{linziqian@fzu.edu.cn}}.}
\and Huangang Lu\thanks{School of Physics, Peking University, Beijing, China.
Email \href{mailto:luhuangang@pku.edu.cn}{\texttt{luhuangang@pku.edu.cn}}.}}
\date{29 October 2026}
\begin{document}
\maketitle

\begin{abstract}
We extend Seginer's matrix norm comparison to tensors
with optimal dependence on the tensor order.
For a real tensor $T$ of order $k\ge2$ with centered, integrable,
independent and identically distributed entries, we prove
\[
 \E f(T)\le\E\norm T_{\inj}
        \le C\sqrt{k\log k}\,\E f(T),
\]
where $f(T)$ is the largest Euclidean norm of a one-dimensional fiber
and $C$ is universal.
The comparison holds for every finite rectangular format, and the
entry law may depend on both the dimensions and the tensor order.
We show that the factor $\sqrt{k\log k}$ is optimal by considering
Rademacher tensors on $[2]^k$.
\end{abstract}



\section{Introduction}\label{sec:intro}
The operator norm of a matrix dominates the Euclidean norm of every
row and column. For an $m\times n$ real random matrix $A$ with
centered integrable i.i.d.\ entries, Seginer's theorem
\cite[Corollary 2.2]{Seginer} gives
\[
 \E\norm A_{\mathrm{op}}
 \asymp \E\max\left\{
      \max_{1\le i\le m}\norm{A_{i,\cdot}}_2,\,
      \max_{1\le j\le n}\norm{A_{\cdot,j}}_2
              \right\},
\]
where $A_{i,\cdot}$ and $A_{\cdot,j}$ denote the rows and columns of
$A$, respectively, and the comparison constants are universal.

In this article, we extend this comparison to tensors of arbitrary
order. For a real tensor $T=(T_{i_1,\ldots,i_k})$ of order $k\ge2$
and format $n_1\times\cdots\times n_k$, where $n_j\ge1$, define
its injective norm by
\[
 \norm T_{\inj}
 =\sup_{\substack{x_j\in\mathbb R^{n_j}\\\norm{x_j}_2\le1}}
   \left|\sum_{i_1=1}^{n_1}\cdots\sum_{i_k=1}^{n_k}
           T_{i_1,\ldots,i_k}x_{1,i_1}\cdots x_{k,i_k}\right|.
\]
Let $f(T)$ be the largest Euclidean norm of a one-dimensional fiber,
obtained by varying one index and fixing the others.
We write $[r]=\{1,\ldots,r\}$ for integers $r\ge1$. All logarithms
are natural. Theorem~\ref{thm:main} gives a comparison uniform over
the dimensions and the entry law.

\begin{theorem}[Uniform fiber comparison]\label{thm:main}
There is a universal constant $C>0$ such that, for every finite
rectangular format and every centered integrable i.i.d.\ real entry law,
\begin{equation}\label{eq:main}
 \E f(T)\le\E\norm T_{\inj}
        \le C\sqrt{k\log k}\,\E f(T).
\end{equation}
The law may depend on $k$ and on all the dimensions.
\end{theorem}

Theorem~\ref{thm:sign} shows that the factor $\sqrt{k\log k}$ is
optimal even for Rademacher tensors with every mode of dimension two.

\begin{theorem}[Sharpness in the tensor order]\label{thm:sign}
Let $R^{(k)}$ have independent Rademacher entries on $[2]^k$.
There are universal $c,C>0$ such that, for every $k\ge2$,
\begin{equation}\label{eq:sign}
 c\sqrt{k\log k}
 \le \E\norm{R^{(k)}}_{\inj}
 \le C\sqrt{k\log k},\qquad f(R^{(k)})=\sqrt2.
\end{equation}
\end{theorem}

The smallest comparison factor at a given tensor order can be written as
\[
 C_k^\ast=
 \sup_T\frac{\E\norm T_{\inj}}{\E f(T)},
\]
where the supremum ranges over all finite rectangular formats of order $k$ and
centered integrable i.i.d.\ entry laws with $\E f(T)>0$.
Theorems~\ref{thm:main} and \ref{thm:sign} give
\begin{equation}\label{eq:optimal}
                C_k^\ast\asymp\sqrt{k\log k}
\end{equation}
with universal comparison constants.

For a scalar random variable $Z$ and $q\ge1$, write
$\norm Z_{L_q}=(\E|Z|^q)^{1/q}$.
Let $M(T)=\max_{i_1,\ldots,i_k}|T_{i_1,\ldots,i_k}|$ be the
largest absolute entry.
Corollary~\ref{cor:higher-moments} extends the expectation bound to
higher moments. Its additional term depends on the largest entry.

\begin{corollary}[Higher moments]\label{cor:higher-moments}
There is a universal constant $C>0$ such that, for every $q\ge1$
and every tensor $T$ as in Theorem~\ref{thm:main} whose entries have
finite $q$th absolute moment,
\begin{equation}\label{eq:moment-refined}
 \bigl\|\norm T_{\inj}\bigr\|_{L_q}
 \le C\left(\sqrt{k\log k}\,\E f(T)
                    +\sqrt q\,\norm{M(T)}_{L_q}\right).
\end{equation}
In particular,
\begin{equation}\label{eq:moment-comparison}
 \norm{f(T)}_{L_q}
 \le \bigl\|\norm T_{\inj}\bigr\|_{L_q}
 \le C\bigl(\sqrt{k\log k}+\sqrt q\bigr)\norm{f(T)}_{L_q}.
\end{equation}
\end{corollary}

The $\sqrt q$ term in \eqref{eq:moment-comparison} cannot be replaced
by $o(\sqrt q)$ uniformly over the tensor order, as shown in
Section~\ref{sec:moments}.

\subsection{Relation to earlier work}
Lucca and Pesenti \cite[Theorem 1.7 and Proposition B.1]{LP}
obtain a fiber comparison for centered i.i.d.\ tensor entries satisfying
\[
 \E|\xi|^{2+\varepsilon}
 \le K(\E\xi^2)^{1+\varepsilon/2},\qquad
 \varepsilon\in(0,1),\quad K\ge1.
\]
The comparison constants depend on $k$, $\varepsilon$, and $K$.
Their Conjecture 1.8 asks whether a comparison depending only on $k$
holds for arbitrary centered integrable i.i.d.\ entries.
To relate their normalization to ours, let $f_j(T)$ be the largest
Euclidean norm of a fiber parallel to mode $j$, and write
$\Phi(T)=\E f(T)$ and $\Psi(T)=\max_j\E f_j(T)$.
Then
\begin{equation}\label{eq:mode-placement}
                \Psi(T)\le\Phi(T)\le k\Psi(T).
\end{equation}
For fixed $k$, these two quantities are therefore equivalent, so
Theorem~\ref{thm:main} establishes their conjecture.
The optimal dependence in \eqref{eq:optimal}
concerns $\Phi(T)$, with the maximum over modes inside the
expectation. The conversion in \eqref{eq:mode-placement} does not
establish the same dependence for $\Psi(T)$.

For independent centered entries with a common subgaussian variance
proxy $\sigma^2$, integrating the tail estimate of Tomioka and Suzuki
\cite[Theorem 1]{TS} gives
\[
 \E\norm T_{\inj}
 \le C\sigma\sqrt{\Bigl(\sum_{j=1}^k n_j\Bigr)\log k}.
\]
This bound depends explicitly on the mode dimensions and the
subgaussian parameter. In Theorem~\ref{thm:main}, the dimensions
enter through $\E f(T)$, while the comparison factor
$C\sqrt{k\log k}$ is independent of all mode dimensions and of
the entry law. The comparison holds for all centered integrable
i.i.d.\ entries.

At a fixed common mode dimension and fixed subgaussian parameter,
the bound of Tomioka and Suzuki already has order $\sqrt{k\log k}$.
Dartois and McKenna obtain finer Gaussian upper bounds using the
Kac--Rice formula \cite{DM24}. Their subsequent moment method
also covers Rademacher entries \cite[Theorem 3.1]{DM26}.
Theorem~\ref{thm:sign} proves a lower bound of the same order on
$[2]^k$, establishing the necessity of this factor in the
dimension-free fiber comparison \eqref{eq:main}.

Complementary results allow entries with different distributions.
Lata\l a, van Handel, and Youssef \cite[Corollary 4.1]{LVY}
prove a fiber comparison for matrices with independent centered
Gaussian entries and arbitrary entry variances.
For independent centered Gaussian or uniformly bounded tensor entries,
Boedihardjo \cite[Theorem 1.1 and Corollary 1.4]{Boe}
bounds the expected injective norm by fiber norms of the
standard-deviation array and an additional logarithmic dimension term.
Proposition~\ref{prop:non-iid} shows that a fiber comparison cannot
hold uniformly over arbitrary independent centered entry laws.
In a different direction, Bandeira et al.\ \cite{BGJLR} and
Aden-Ali \cite{AA} study tensor series with Gaussian and subgaussian
coefficients, allowing correlations among entries and using
variance parameters to control the injective norm.

For sparse Bernoulli tensors, Zhou and Zhu \cite{ZZ} establish
matching expectation bounds, concentration, and regularization
results at fixed tensor order.
The application in Section~\ref{sec:applications} recovers the
expected norm scale at parameter $1/n$ from the general comparison
and elementary estimates for the largest fiber.

\subsection{Proof strategy and organization}
Symmetrization and truncation reduce the proof of the upper bound
in Theorem~\ref{thm:main} to the case of bounded symmetric entries.
The truncation level is chosen so that the entries exceeding it
in absolute value form a tensor with expected injective norm at
most $C\E f(T)$. For the bounded truncation, a product net controls
small summands in multilinear contractions. It remains to control
the large summands through the nonnegative process in
\eqref{eq:soft}.

The main step is a finite covering of the event that this process
is large while the fiber norm is bounded. The covering events are
indexed by a deterministic finite family of product vectors, and
each is a tail event for a sum of independent bounded terms.
The fiber bound enters the construction through selection and
rounding, but is absent from the covering events.
Bennett's inequality bounds their probabilities, and the
cardinality estimate controls the union bound. Integrating the
resulting tail estimate gives the upper bound with factor
$\sqrt{k\log k}$.

For the lower bound in Theorem~\ref{thm:sign}, we construct a
Euclidean packing of unit product vectors with separation bounded
below by a positive universal constant and logarithmic cardinality
of order $k\log k$.
Sudakov minoration gives the corresponding Gaussian lower bound
of order $\sqrt{k\log k}$.
We choose the product vectors to have exponentially small
coordinates so that replacing Gaussian entries by Rademacher
signs introduces an error negligible relative to this lower bound.

Section~\ref{sec:reduction} derives Theorem~\ref{thm:main} from
Theorem~\ref{thm:prob}, which is proved in
Section~\ref{sec:large-contractions}.
Section~\ref{sec:sharpness} proves Theorem~\ref{thm:sign}, and
Section~\ref{sec:applications} contains the higher-moment bound,
the Bernoulli application, and a matrix counterexample with
non-identically distributed entries.
\section{Reduction to a nonnegative process}\label{sec:reduction}
This section derives Theorem~\ref{thm:main} from a probability
estimate for the large summands in multilinear contractions.
The estimate is stated in Theorem~\ref{thm:prob} and proved in
Section~\ref{sec:large-contractions}.
Symmetrization and truncation allow us to apply it under the
first-moment assumption of Theorem~\ref{thm:main}.

Write $n=\max_j n_j$ for the maximum fiber length and
$d=\sum_{j=1}^k n_j$ for the total number of mode-labelled coordinates.
We index entries by multi-indices
$e=(e_1,\ldots,e_k)\in\mathcal I=\prod_{j=1}^k[n_j]$, so that
$T_e=T_{e_1,\ldots,e_k}$.

We identify each tensor with its associated multilinear form.
For a tuple of factors $x=(x_1,\ldots,x_k)$, use $x_e$ for its
product coordinate at $e$, so that
\[
 T(x_1,\ldots,x_k)=\sum_{e\in\mathcal I}T_ex_e,\qquad
 x_e=\prod_jx_{j,e_j}.
\]
Thus $x_e$ is the $e$th coordinate of $x_1\otimes\cdots\otimes x_k$.
The large summands will be measured by their excess above a threshold.
For a nonnegative tensor $A=(A_e)_{e\in\mathcal I}$ and $t>0$, define
\begin{equation}\label{eq:soft}
 H_t(A)=\max_{\substack{x_j\ge0\\\norm{x_j}_2\le1}}
                 \sum_e(A_ex_e-t)_+.
\end{equation}
Here $(s)_+=\max\{s,0\}$, so only summands exceeding $t$ contribute.
Continuity and compactness ensure that the maximum is attained.
Theorem~\ref{thm:prob} controls this process on the event that the
largest fiber norm is bounded.

\begin{theorem}[Nonnegative process estimate]\label{thm:prob}
There are universal constants $c,C>0$ with the following property.
Let $A=(A_e)_{e\in\mathcal I}$ have independent nonnegative entries
satisfying $\E A_e^2\le1/n$, and set $t=1/(n\sqrt{k\log k})$.
For every $F\ge1$ and $u\ge CF\sqrt{k\log k}$,
\begin{equation}\label{eq:tail}
 \PP\{H_t(A)>u,\ f(A)\le F\}
 \le \exp\left\{-c\frac{u^2}{kF^2}\log d\right\}.
\end{equation}
Consequently,
\begin{equation}\label{eq:soft-expectation}
 \E H_t(A)\le C\sqrt{k\log k}\{1+\E[f(A)^2]\}.
\end{equation}
\end{theorem}

\subsection{A norm bound for symmetric entries}
Lemma~\ref{lem:soft-norm} relates the expected injective norm of a
tensor with symmetric entries to the nonnegative process $H_t$
applied to its entrywise absolute values. This allows us to use
Theorem~\ref{thm:prob} to estimate the tensor norm.

\begin{lemma}\label{lem:soft-norm}
Let $W=(W_e)_{e\in\mathcal I}$ be a real random tensor with independent
entries whose distributions are symmetric about zero.
Suppose that $\E W_e^2\le1/n$, and set $t=1/(n\sqrt{k\log k})$.
Then, for a universal constant $C>0$,
\begin{equation}\label{eq:soft-norm}
                 \E\norm W_{\inj}\le C\sqrt{k\log k}+4\E H_t(|W|).
\end{equation}
If also $\E[f(W)^2]\le1$, then
\begin{equation}\label{eq:normalized}
                         \E\norm W_{\inj}\le C\sqrt{k\log k}.
\end{equation}
\end{lemma}
\begin{proof}
The Euclidean covering estimate in \cite[Corollary 4.2.13]{HDP}
gives a $1/(2k)$-net $\mathcal N_j$ of the unit sphere in
$\mathbb R^{n_j}$ with $|\mathcal N_j|\le(1+4k)^{n_j}$.
Put $\mathcal N=\mathcal N_1\times\cdots\times\mathcal N_k$ and
$M=|\mathcal N|$. Then $M\le(1+4k)^d$ and
$\log(2M)\le4nk\log k$.
For any unit factors $x_1,\ldots,x_k$, choose corresponding net points
$y_j\in\mathcal N_j$ with $\norm{x_j-y_j}_2\le1/(2k)$.
Replacing the factors one at a time and using multilinearity gives
\[
 |W(x_1,\ldots,x_k)-W(y_1,\ldots,y_k)|
 \le\norm W_{\inj}\sum_j\norm{x_j-y_j}_2
 \le\frac12\norm W_{\inj}.
\]
By the triangle inequality and then taking the supremum over all
$x_j\in\mathbb R^{n_j}$ with $\norm{x_j}_2=1$, we obtain
\[
 \norm W_{\inj}
 \le\max_{y\in\mathcal N}|W(y_1,\ldots,y_k)|+\frac12\norm W_{\inj}.
\]
Moving the last term to the left gives
\begin{equation}\label{eq:net-approximation}
 \norm W_{\inj}\le2\max_{y\in\mathcal N}|W(y_1,\ldots,y_k)|.
\end{equation}
It therefore suffices to control contractions on this finite set.
For $x=(x_1,\ldots,x_k)\in\mathcal N$, write
$W(x)=W(x_1,\ldots,x_k)$ and define the truncated sum
\[
 W_{\le2t}(x)=\sum_{e\in\mathcal I}
          W_ex_e\one_{\{|W_ex_e|\le2t\}}.
\]
For fixed $x$, the summands are independent because the summand indexed
by $e$ depends only on $W_e$.
The map $u\mapsto u\one_{\{|u|\le2t\}}$ is odd, and $W_e$ has a
distribution symmetric about zero. Thus each truncated summand has mean zero,
\[
 \E\!\left[W_ex_e\one_{\{|W_ex_e|\le2t\}}\right]=0.
\]
Since every $x_j$ is a unit vector,
$\sum_e x_e^2=\prod_j\norm{x_j}_2^2=1$.
Independence and centering therefore give
\begin{align*}
 \Var(W_{\le2t}(x))
 &=\sum_{e\in\mathcal I}
       \E\!\left[W_e^2x_e^2\one_{\{|W_ex_e|\le2t\}}\right]\\
 &\le\sum_{e\in\mathcal I}x_e^2\E W_e^2
 \le\frac1n\sum_{e\in\mathcal I}x_e^2=\frac1n.
\end{align*}
Since the summands are bounded in absolute value by $2t$,
Bernstein's exponential moment bound \cite[Theorem~2.10]{BLM},
applied to $W_{\le2t}(x)$ and its negative, gives
\[
 \log\E e^{\lambda W_{\le2t}(x)}
 \le\frac{\lambda^2}{2n(1-2t|\lambda|/3)},
 \qquad |\lambda|<\frac{3}{2t}.
\]
Since $e^{\lambda\max_{x\in\mathcal N}|W_{\le2t}(x)|}
\le\sum_{x\in\mathcal N}(e^{\lambda W_{\le2t}(x)}+e^{-\lambda W_{\le2t}(x)})$,
Jensen's inequality and the preceding estimate give, for $0<\lambda<3/(2t)$,
\begin{align*}
 \E\max_{x\in\mathcal N}|W_{\le2t}(x)|
 &\le\frac1\lambda\log\E e^{\lambda\max_{x\in\mathcal N}|W_{\le2t}(x)|}\\
 &\le\frac1\lambda\log\sum_{x\in\mathcal N}
               \bigl(\E e^{\lambda W_{\le2t}(x)}+\E e^{-\lambda W_{\le2t}(x)}\bigr)
 \le\frac{\log(2M)}{\lambda}
       +\frac{\lambda}{2n(1-2t\lambda/3)}.
\end{align*}
Optimizing over $\lambda$ gives
$\lambda^\ast=((2n\log(2M))^{-1/2}+2t/3)^{-1}$, and hence
\begin{align*}
 \E\max_{x\in\mathcal N}|W_{\le2t}(x)|
 &\le\sqrt{2\log(2M)/n}+(2t/3)\log(2M)\\
 &\le6\sqrt{k\log k}.
\end{align*}
The remaining terms satisfy $|W_ex_e|>2t$.
Using $w\one_{\{w>2t\}}\le2(w-t)_+$ for $w\ge0$, we obtain,
for every $x\in\mathcal N$,
\begin{align*}
 |W(x)-W_{\le2t}(x)|
 &\le\sum_e|W_e|\,|x_e|\one_{\{|W_ex_e|>2t\}}\\
 &\le2\sum_e\bigl(|W_e|\,|x_e|-t\bigr)_+\\
 &\le2H_t(|W|).
\end{align*}
The last inequality follows by taking the coordinatewise absolute
values $|x_j|$ as the nonnegative unit factors in \eqref{eq:soft}.
Together with \eqref{eq:net-approximation}, this yields
\begin{align*}
 \E\norm W_{\inj}
 &\le2\E\max_{x\in\mathcal N}|W_{\le2t}(x)|+4\E H_t(|W|)\\
 &\le C\sqrt{k\log k}+4\E H_t(|W|),
\end{align*}
which proves \eqref{eq:soft-norm}.
When $\E[f(W)^2]\le1$, Theorem~\ref{thm:prob} gives
$\E H_t(|W|)\le C\sqrt{k\log k}$.
Substituting this into \eqref{eq:soft-norm} proves \eqref{eq:normalized}
after increasing $C$ if necessary.
\end{proof}

\subsection{Symmetrization and truncation}
To apply Lemma~\ref{lem:soft-norm} to a bounded truncation, we control
its fiber second moment using the variance estimate in
Lemma~\ref{lem:fiber-var}.

\begin{lemma}\label{lem:fiber-var}
Let $Y=(Y_e)_{e\in\mathcal I}$ be a real random tensor with independent
entries satisfying $|Y_e|\le b$ almost surely. Then
\[
                            \Var f(Y)\le b^2.
\]
\end{lemma}
\begin{proof}
Let $Y^{(e)}$ be obtained from $Y$ by replacing $Y_e$ with an independent copy.
The Efron--Stein inequality \cite[Proposition 1]{BBLM} and
exchangeability give
\[
 \Var f(Y)\le\frac12\sum_e\E\{f(Y)-f(Y^{(e)})\}^2
            =\sum_e\E\{f(Y)-f(Y^{(e)})\}_+^2.
\]
Fix an ordering of the fibers and let $J\subseteq\mathcal I$ be the
index set of the first fiber whose Euclidean norm is $f(Y)$.
Using the same index set $J$ for $Y^{(e)}$, we have
\[
 f(Y)^2-Y_e^2\one_{\{e\in J\}}
 \le f(Y)^2+\bigl((Y_e^{(e)})^2-Y_e^2\bigr)\one_{\{e\in J\}}
 =\sum_{a\in J}(Y_a^{(e)})^2
 \le f(Y^{(e)})^2.
\]
On the event $\{f(Y)>0\}$, factorization gives
\[
 \{f(Y)-f(Y^{(e)})\}_+
 =\frac{\{f(Y)^2-f(Y^{(e)})^2\}_+}{f(Y)+f(Y^{(e)})}
 \le\frac{\{f(Y)^2-f(Y^{(e)})^2\}_+}{f(Y)}
 \le\frac{Y_e^2}{f(Y)}\one_{\{e\in J\}},
\]
where the first inequality uses $f(Y^{(e)})\ge0$.
Squaring and summing yields
\[
 \sum_e\{f(Y)-f(Y^{(e)})\}_+^2
 \le\frac{\sum_{e\in J}Y_e^4}{f(Y)^2}
 \le\frac{b^2\sum_{e\in J}Y_e^2}{f(Y)^2}=b^2.
\]
The last equality uses $\sum_{e\in J}Y_e^2=f(Y)^2$, by the choice of $J$.
When $f(Y)=0$, all the positive parts vanish.
Taking expectations completes the proof.
\end{proof}

\begin{proof}[Proof of Theorem~\ref{thm:main}]
For a fiber vector $v$ in mode $j$, take the other factors to be the
standard basis vectors corresponding to its fixed indices. Then
\[
 \norm v_2=\sup_{\norm{x_j}_2\le1}|\langle v,x_j\rangle|
 \le\norm T_{\inj}.
\]
Maximizing over the fibers gives $f(T)\le\norm T_{\inj}$.
Both norms are integrable because
$\norm T_{\inj}\le\sum_e|T_e|$.

Let $T'$ be an independent copy of $T$ and put $S=T-T'$.
The entries of $S$ are symmetric i.i.d., and $\E[S\mid T]=T$.
Jensen's inequality and the triangle inequality give
\begin{equation}\label{eq:symmetrization}
 \E\norm T_{\inj}\le\E\norm S_{\inj},\qquad
 \E f(S)\le2\E f(T).
\end{equation}
If $\E f(S)=0$, then $S=0$ almost surely and $\E\norm T_{\inj}=0$,
so the upper bound holds. Otherwise, divide $T$, $T'$ and $S$ by
$\E f(S)$, retaining the same notation. This preserves $S=T-T'$
and \eqref{eq:symmetrization}. Since the desired upper bound is
homogeneous, we may therefore assume $\E f(S)=1$.
Under this normalization, define
\[
             Y_e=S_e\one_{\{|S_e|\le2\}},\qquad U=S-Y.
\]
Put $N=|\mathcal I|$ and $p=\PP(|S_e|>2)$. Markov's inequality gives
\[
 \PP\!\left(\max_e|S_e|>2\right)
 \le\frac{\E\max_e|S_e|}{2}
 \le\frac{\E f(S)}2=\frac12.
\]
Hence $\frac12\le\PP(\max_e|S_e|\le2)=(1-p)^N
\le(1-p)^{N-1}$.
Let $B$ be the event that exactly one entry of $U$ is nonzero.
Since $\PP(U_e=0)=1-p$, independence gives
\[
 \E\!\left[\max_e|U_e|\one_B\right]
 =\sum_e\E\!\left[|U_e|\prod_{a\ne e}\one_{\{U_a=0\}}\right]
 =N\E|U_e|(1-p)^{N-1}.
\]
It follows that
\[
 \frac N2\E|U_e|\le N\E|U_e|(1-p)^{N-1}
 \le\E\max_e|U_e|\le\E f(S)=1.
\]
Hence
\begin{equation}\label{eq:tail-transfer}
                   \E\norm U_{\inj}\le\sum_e\E|U_e|\le2.
\end{equation}

The entries of $Y$ are symmetric i.i.d.
Using a fiber of length $n$, Lemma~\ref{lem:fiber-var}, and
$\E f(Y)\le1$, we obtain
\[
 n\E Y_e^2\le\E[f(Y)^2]
 =\Var f(Y)+(\E f(Y))^2\le5.
\]
Thus $W=Y/\sqrt5$ satisfies both
$n\E W_e^2\le1$ and $\E[f(W)^2]\le1$.
Applying Lemma~\ref{lem:soft-norm} to $W$ and using
\eqref{eq:tail-transfer} gives
\[
 \E\norm S_{\inj}
 \le\E\norm Y_{\inj}+\E\norm U_{\inj}
 \le C\sqrt{k\log k}+2
 \le C\sqrt{k\log k}.
\]
Combining this estimate with \eqref{eq:symmetrization} and
$\E f(S)=1$ yields
\[
 \E\norm T_{\inj}\le\E\norm S_{\inj}
 \le C\sqrt{k\log k}\,\E f(S)
 \le C\sqrt{k\log k}\,\E f(T).
\]
Rescaling $T$ proves the upper bound in general.
\end{proof}

\section{A finite cover for large contractions}\label{sec:large-contractions}
We first construct a deterministic finite cover of the event in
Theorem~\ref{thm:prob}, then estimate the probabilities of the
finitely many covering events.

\subsection{Construction of the finite cover}\label{sec:covering}
Fix $u,F,t>0$. The construction proceeds in three steps.
First, we use a maximizing tuple for $H_t(A)$ to compute coordinate
contributions and derive weighted estimates. Second, we select
entries and coordinates to preserve a large weighted sum while
bounding individual terms and keeping coordinates large enough
for rounding. Third, we round the surviving coordinates to obtain
a finite cover of $\{H_t(A)>u,\ f(A)\le F\}$.

\subsubsection{Maximizing directions and weighted contributions}\label{sec:weighted-contributions}
For nonnegative tensors $A$, the definition of $H_t(A)$ gives
\begin{equation}\label{eq:direction-decomposition}
 \{H_t(A)>u,\ f(A)\le F\}
 =\bigcup_{\substack{x_j\ge0,\ \norm{x_j}_2=1\\1\le j\le k}}
 \left\{\sum_e(A_ex_e-t)_+>u,\ f(A)\le F\right\}.
\end{equation}
Indeed, for $A\in\{H_t(A)>u,\ f(A)\le F\}$, there are nonnegative
factors with $\norm{x_j}_2\le1$ and
$\sum_e(A_ex_e-t)_+>u$. Since $u>0$, every $x_j$ is nonzero.
Set $y_j=x_j/\norm{x_j}_2$, so that $\norm{y_j}_2=1$. Then
\[
 y_e=\frac{x_e}{\prod_j\norm{x_j}_2}\ge x_e,\qquad
 \sum_e(A_ey_e-t)_+\ge\sum_e(A_ex_e-t)_+>u.
\]
The fiber condition is unchanged, so $A$ belongs to the union on
the right. Conversely, if $A$ belongs to that union, some nonnegative
unit tuple $x$ satisfies $f(A)\le F$ and
\[
 H_t(A)\ge\sum_e(A_ex_e-t)_+>u.
\]
This gives the inclusion from right to left and proves
\eqref{eq:direction-decomposition}.

The union in \eqref{eq:direction-decomposition} is uncountable whenever
some $n_j\ge2$. To obtain a finite cover, we will select and round
coordinates. Setting $x_{j,i}=0$ removes the contribution
\[
 \sum_{e:\,e_j=i}(A_ex_e-t)_+
 \le\sum_{\substack{e:\,A_ex_e>t\\e_j=i}}A_ex_e.
\]
Lemma~\ref{lem:coordinate-contributions} evaluates the right-hand side
at a maximizing tuple.

\begin{lemma}[Coordinate contributions at a maximizer]\label{lem:coordinate-contributions}
Let $A\ge0$ and $t>0$, with $H_t(A)>0$. Suppose that
$x=(x_1,\ldots,x_k)$ is a tuple of nonnegative unit factors satisfying
\[
 \sum_e(A_ex_e-t)_+=H_t(A).
\]
Set
\[
 E=\{e:A_ex_e>t\},\qquad
 S=\sum_{e\in E}A_ex_e=H_t(A)+t|E|.
\]
Then, for every tuple $y$ of nonnegative unit factors,
\[
 \sum_{e\in E}A_ey_e\le S=\sum_{e\in E}A_ex_e.
\]
Moreover, for every $1\le j\le k$ and $i\in[n_j]$,
\begin{equation}\label{eq:coordinate-contribution}
 \sum_{\substack{e\in E\\e_j=i}}A_ex_e=Sx_{j,i}^2.
\end{equation}
\end{lemma}
\begin{proof}
Keep the set $E$ defined by $x$ fixed. For every nonnegative unit
tuple $y$,
\[
 \sum_{e\in E}A_ey_e
 \le\sum_e(A_ey_e-t)_++t|E|
 \le H_t(A)+t|E|=S.
\]
Equality holds at $y=x$, so $x$ also maximizes the multilinear
form $\sum_{e\in E}A_ey_e$.

Fix $j$ and keep all other factors unchanged. As a function of
$y_j$, this form is $\langle g,y_j\rangle$, where
\[
 g_i=\sum_{\substack{e\in E\\e_j=i}}
       A_e\prod_{\ell\ne j}x_{\ell,e_\ell}.
\]
Since $g\ge0$ and $\langle g,x_j\rangle=S>0$, the vector
$g/\norm g_2$ is an admissible replacement for $x_j$. Therefore
\[
 S=\langle g,x_j\rangle
 \le\norm g_2
 =\left\langle g,\frac{g}{\norm g_2}\right\rangle
 \le S.
\]
The first inequality is Cauchy--Schwarz, and the last follows from
the maximality of $x$. Equality throughout gives $g=Sx_j$.
Multiplying its $i$th component by $x_{j,i}$ proves
\eqref{eq:coordinate-contribution}.
\end{proof}

Fix $A\in\{H_t(A)>u,\ f(A)\le F\}$. By continuity and compactness,
the maximum over the unit tuples in \eqref{eq:direction-decomposition}
is attained. Choose a maximizing tuple $x$, and let $E$ and $S>u$
be as in Lemma~\ref{lem:coordinate-contributions}.
Dividing \eqref{eq:coordinate-contribution} by $x_{j,i}^2$ makes
each coordinate used by $E$ contribute $S$. This leads to the weights
\[
 w_e=\sum_{j=1}^k\frac1{x_{j,e_j}^2}\quad(e\in E).
\]
For $K\subseteq E$, let $\pi_jK=\{e_j:e\in K\}\subseteq[n_j]$
be its projection onto the $j$th mode, recording the coordinates
of $x_j$ used by the entries in $K$.
Summing the coordinate identities gives
\begin{equation}\label{eq:weighted-total}
 \sum_{e\in E}A_ex_ew_e
 =\sum_{j=1}^k\sum_{i\in\pi_jE}\frac1{x_{j,i}^2}
       \sum_{\substack{e\in E\\e_j=i}}A_ex_e
 =S\sum_{j=1}^k|\pi_jE|.
\end{equation}
Each reciprocal $x_{j,e_j}^{-2}$ cancels one coordinate square, so
the fiber condition also gives
\begin{equation}\label{eq:fiber-energy}
 \sum_{e\in E}(A_ex_e)^2w_e
 =\sum_j\sum_{e\in E}A_e^2\prod_{\ell\ne j}x_{\ell,e_\ell}^2
 \le kF^2.
\end{equation}
The estimates \eqref{eq:weighted-total} and \eqref{eq:fiber-energy}
will be used in Section~\ref{sec:entry-coordinate-selection}.

\subsubsection{Selection of entries and coordinates}\label{sec:entry-coordinate-selection}
To construct the retained set, we first discard entries with large
$A_ex_ew_e$ and then select coordinates. We measure the loss at each
stage as a fraction of the total contraction $S$. Accordingly, define
\[
 \mu_e=\frac{A_ex_e}{S}\quad(e\in E),\qquad
 \sum_{e\in E}\mu_e=1.
\]
Thus $\mu_e$ is the fraction of $S$ contributed by entry $e$.
Summing over the entries that use coordinate $(j,i)$ gives
\[
 p_{j,i}=\sum_{\substack{e\in E\\e_j=i}}\mu_e
 =\frac1S\sum_{\substack{e\in E\\e_j=i}}A_ex_e
 =x_{j,i}^2.
\]
The last equality is \eqref{eq:coordinate-contribution}.
The measure $\mu$ is a probability measure on $E$, and
$p_j=(p_{j,i})_{i\in\pi_jE}$ is its $j$th coordinate marginal.

For a retained set $K\subseteq E$, continue to divide contributions
by the original total $S$ and put
\[
 \mu(K)=\sum_{e\in K}\mu_e,\qquad
 p_{j,i}(K)=\sum_{\substack{e\in K\\e_j=i}}\mu_e.
\]
Here $\mu(K)$ is the fraction of $S$ retained in total, and
$p_{j,i}(K)$ is the part retained through coordinate $(j,i)$.
In particular,
\[
 \frac{p_{j,i}(K)}{p_{j,i}}
 =\frac{\displaystyle\sum_{\substack{e\in K\\e_j=i}}A_ex_e}
          {Sx_{j,i}^2}\quad(i\in\pi_jK)
\]
is the fraction of the original contribution through coordinate
$(j,i)$ that remains.

We impose two conditions on the retained coordinates. A lower
bound on $p_{j,i}=x_{j,i}^2$ will control the relative error when
these squared coordinates are rounded. To preserve the weighted
sum, it is enough that each coordinate retain a fixed fraction
$0<c_1<1$ of its original contribution. Indeed, if
$p_{j,i}(K)\ge c_1 p_{j,i}$ for every $i\in\pi_jK$, then
\[
 \sum_{e\in K}A_ex_ew_e
 =S\sum_{j=1}^k\sum_{i\in\pi_jK}\frac{p_{j,i}(K)}{p_{j,i}}
 \ge c_1S\sum_{j=1}^k|\pi_jK|.
\]
It remains to construct a set satisfying these coordinate
conditions and an upper bound on every weighted term. We prepare
the estimates for the two deletion stages.
For entry deletion, normalize \eqref{eq:fiber-energy} to obtain
\[
 \mathcal E=\sum_{e\in E}\mu_e^2w_e
 =\frac1{S^2}\sum_{e\in E}(A_ex_e)^2w_e
 \le\frac{kF^2}{S^2}.
\]
Removing entries with $\mu_ew_e>4\mathcal E$ loses at most $1/4$
of the mass by Markov's inequality.

Deleting a coordinate $(j,i)$ with $p_{j,i}(K)<c_1 p_{j,i}$
decreases $\mu(K)$ by $p_{j,i}(K)$.
Lemma~\ref{lem:coordinate-contributions} gives
$\sum_{e\in E}A_ey_e\le S$ for every tuple of nonnegative unit factors.
For $X=(X_1,\ldots,X_k)$ with law $\mu$ and nonnegative functions
$b_j$ on $\pi_jE$, write
$x_jb_j=(x_{j,i}b_j(i))_{i\in\pi_jE}$, with zeros on the remaining
coordinates. Apply this inequality to the unit factors
$y_j=x_jb_j/\norm{x_jb_j}_2$, provided all these norms are positive.
This gives
\[
 \E_\mu\prod_j b_j(X_j)
 =\frac1S\sum_{e\in E}A_e\prod_j
       \bigl(x_{j,e_j}b_j(e_j)\bigr)
 \le\prod_j\norm{x_jb_j}_2
 =\prod_j(\E_{p_j}b_j^2)^{1/2}.
\]
If any factor has norm zero, both sides vanish.

The coordinates used by $K$ form the set
\[
 V(K)=\{(j,i):1\le j\le k,\ i\in\pi_jK\},\qquad
 |V(K)|=\sum_j|\pi_jK|.
\]
In terms of the marginals, $w_e=\sum_jp_{j,e_j}^{-1}$.
Lemma~\ref{lem:core} constructs the retained set needed for rounding.

\begin{lemma}[Selection of entries and coordinates]\label{lem:core}
Let $\mu$ be a probability measure with positive weights on a
nonempty set $E\subseteq\prod_{j=1}^k[n_j]$ with $|V(E)|\le d$.
Suppose that, for $X=(X_1,\ldots,X_k)$ with law $\mu$ and all
nonnegative functions $b_j$ on $\pi_jE$,
\begin{equation}\label{eq:product}
 \E_\mu\prod_j b_j(X_j)
 \le\prod_j(\E_{p_j}b_j^2)^{1/2}.
\end{equation}
For sufficiently small universal constants $c_1,c_2>0$, there is a
nonempty $J\subseteq E$, with
$v=|V(J)|=\sum_{j=1}^k|\pi_jJ|$, such that
\begin{align}
 \mu(J)&>\frac12, \label{eq:mass}\\
 p_{j,i}(J)&\ge c_1 p_{j,i}\quad((j,i)\in V(J)), \label{eq:retention}\\
 p_{j,i}&\ge\frac{c_2}{v\log^2(ed/v)}\quad((j,i)\in V(J)), \label{eq:floor}\\
 \mu_ew_e&\le4\mathcal E\quad(e\in J). \label{eq:score}
\end{align}
The quantities $p_{j,i}$ and $\mathcal E$ are computed from $\mu$ on $E$.
The restriction to $J$ is not normalized.
\end{lemma}
\begin{proof}
The proof has three steps. First, we select entries to obtain
\eqref{eq:score}. Second, we delete coordinates until
\eqref{eq:retention} and \eqref{eq:floor} hold. Third, we bound the
total deleted mass to prove \eqref{eq:mass} and ensure that the
terminal set is nonempty.
Let $0<c_1,c_2<1/2$ be chosen at the end of the proof.

\emph{Step 1. Entry selection.} Set
$J_0=\{e\in E:\mu_ew_e\le4\mathcal E\}$. Markov's inequality gives
\[
 1-\mu(J_0)
  \le\frac1{4\mathcal E}\sum_{e\in E}\mu_e^2w_e=\frac14.
\]
Every remaining entry satisfies \eqref{eq:score}, and
$\mu(J_0)\ge3/4$.

\emph{Step 2. Coordinate selection.} Start with $K=J_0$ and let
$r=|V(K)|=\sum_j|\pi_jK|$ be the current number of coordinates.
Deleting $(j,i)$ decreases $\mu(K)$ by $p_{j,i}(K)\le p_{j,i}$
and reduces $r$.
We use the size cutoff $c_2/[r\log^2(ed/r)]$, whose sum over
$1\le r\le d$ is less than $5c_2$ by \eqref{eq:summable}.

At each step, delete a coordinate $(j,i)\in V(K)$ if
\[
 \text{(i)}\quad p_{j,i}<\frac{c_2}{r\log^2(ed/r)}
 \qquad\text{or}\qquad
 \text{(ii)}\quad p_{j,i}(K)<c_1 p_{j,i}.
\]
Rule (i) removes coordinates that are too small for rounding.
Rule (ii) removes coordinates retaining less than a fraction $c_1$ of their
original contribution. Use rule (i) whenever possible.
Deleting $(j,i)$ means removing all current entries that use it,
\[
 K\longleftarrow K\setminus\{e\in K:e_j=i\}.
\]
Repeat the two tests with the updated $K$ and $r=|V(K)|$,
stopping if $K=\varnothing$ or neither rule applies.
Each update reduces $r$, so the procedure terminates and no
coordinate is chosen twice.
Let $J$ be the terminal set.

\emph{Step 3. Retained mass and nonemptiness.} We prove that
$\mu(K)>1/2$ throughout the procedure. Suppose, to the contrary,
that $K$ is the first remaining set with $\mu(K)\le1/2$.
Let $K^{-}$ be the set immediately before the deletion producing
$K$, and let $(j,i)$ be the deleted coordinate. By the choice of $K$,
\[
 \mu(K^{-})>\frac12,\qquad
 \mu(K^{-})-\mu(K)
 =\sum_{\substack{e\in K^{-}\\e_j=i}}\mu_e
 =p_{j,i}(K^{-}).
\]
Under rule (i), $p_{j,i}(K^{-})\le p_{j,i}<c_2$, since
$r\log^2(ed/r)\ge1$ for $1\le r\le d$.
Under rule (ii), $p_{j,i}(K^{-})<c_1p_{j,i}\le c_1$.
Consequently,
\[
 \mu(K)=\mu(K^{-})-p_{j,i}(K^{-})
 >\frac12-\max\{c_1,c_2\}>0.
\]
In particular, $K$ is nonempty.

We first bound the total mass removed under rule (ii). Every
coordinate $(j,i)$ selected under this rule is absent from $V(K)$,
and its deletion removes mass less than $c_1 p_{j,i}$.
Since no coordinate is selected twice, the total mass removed is at most
$c_1\sum_{(j,i)\in V(E)\setminus V(K)}p_{j,i}$.
To bound this sum, apply \eqref{eq:product} with
$b_j=\one_{\pi_jK}$ to obtain
\[
 \mu(K)\le\E_\mu\prod_j\one_{\pi_jK}(X_j)
 \le\prod_jp_j(\pi_jK)^{1/2}.
\]
Since $1-s\le-\log s$ for $s>0$,
\[
 \sum_{(j,i)\in V(E)\setminus V(K)}p_{j,i}
 =\sum_j\{1-p_j(\pi_jK)\}
 \le-\sum_j\log p_j(\pi_jK)
 \le2\log\frac1{\mu(K)}.
\]
Thus the total mass removed under rule (ii) is at most
$2c_1\log(1/\mu(K))$.

For rule (i), a deletion with $r$ current coordinates removes mass
less than $c_2/[r\log^2(ed/r)]$. Since $r$ strictly decreases,
summing these thresholds over $1\le r\le d$ bounds the total
mass removed under rule (i).
Grouping $r$ into the intervals $(d2^{-m-1},d2^{-m}]$, $m\ge0$,
gives
\begin{equation}\label{eq:summable}
 \sum_{r=1}^d\frac1{r\log^2(ed/r)}
 \le2\sum_{m\ge0}\frac1{(1+m\log2)^2}
 \le2\left(1+\frac1{\log2}\right)<5.
\end{equation}
Thus rule (i) removes mass less than $5c_2$.

Since $\mu(J_0)\ge3/4$ and $\mu(K)\le1/2$, combining the two
bounds on the removed mass gives
\[
\begin{aligned}
 \frac14\le\mu(J_0)-\mu(K)
 &\le c_2\sum_{r=1}^d\frac1{r\log^2(ed/r)}
      +c_1\sum_{(j,i)\in V(E)\setminus V(K)}p_{j,i}\\
 &<5c_2+2c_1\log\frac1{\mu(K)}.
\end{aligned}
\]
We now fix universal constants $c_1,c_2>0$ satisfying
\[
 \max\{c_1,c_2\}\le\frac14,\qquad 5c_2+2c_1\log4<\frac14.
\]
The bound for the last deletion now gives
$\mu(K)>1/2-\max\{c_1,c_2\}\ge1/4$. Substituting this into the
estimate for the removed mass yields
\[
 \frac14\le\mu(J_0)-\mu(K)
 <5c_2+2c_1\log\frac1{\mu(K)}
 <5c_2+2c_1\log4<\frac14,
\]
which is impossible.

Thus the mass remains greater than $1/2$ throughout the procedure.
In particular, $\mu(J)>1/2$, proving \eqref{eq:mass} and giving
$J\ne\varnothing$. The procedure therefore stops because neither
deletion rule applies. With $v=|V(J)|$, the stopping conditions give
\eqref{eq:retention} and \eqref{eq:floor}, while $J\subseteq J_0$
gives \eqref{eq:score}. This proves Lemma~\ref{lem:core}.
\end{proof}

We now use Lemma~\ref{lem:core} to bound the retained weighted sum
and its individual terms, and to prepare the coordinates for rounding.
Apply the lemma to $\mu_e=A_ex_e/S$, whose product inequality
\eqref{eq:product} was verified before the lemma. Let $J$ be the
resulting subset and $v=|V(J)|=\sum_j|\pi_jJ|$.

First, $\mu_e=A_ex_e/S$ and $p_{j,i}=x_{j,i}^2$ express the
weighted sum through the retained marginals. Using
\eqref{eq:retention} and $S>u$ gives
\begin{equation}\label{eq:core-sum}
 \sum_{e\in J}A_ex_ew_e
 =S\sum_{j=1}^k\sum_{i\in\pi_jJ}\frac{p_{j,i}(J)}{p_{j,i}}
 \ge c_1Sv>c_1uv.
\end{equation}
For the individual terms, \eqref{eq:score} and
\eqref{eq:fiber-energy} give
\begin{equation}\label{eq:core-bound}
 A_ex_ew_e=S\mu_ew_e
 \le4S\mathcal E\le\frac{4kF^2}{S}<\frac{4kF^2}{u}
 \qquad(e\in J).
\end{equation}
Finally, to prepare for rounding, set $I_j=\pi_jJ$ for each mode
$j$. Since $J$ is nonempty, every $I_j$ is nonempty. The bounds
\eqref{eq:mass} and \eqref{eq:floor} give
\[
 p_j(I_j)\ge\mu(J)>\frac12,\qquad
 p_{j,i}\ge\frac{c_2}{v\log^2(ed/v)}\quad(i\in I_j).
\]
The bound $p_j(I_j)>1/2$ controls the normalization within each
retained factor. The lower bound on $p_{j,i}$ controls the relative
error when these squared coordinates are rounded.

\subsubsection{Rounding and finite covering}\label{sec:rounding-cover}
We discretize the retained coordinate weights $p_{j,i}=x_{j,i}^2$,
which may take infinitely many values as $A$ varies.
Lemma~\ref{lem:covering} rounds and normalizes these weights to
$z_{j,i}^2$ on the same supports
$\operatorname{supp}(z_j)=I_j=\pi_jJ$.
The resulting tuples $z=(z_1,\ldots,z_k)$ belong to finite
families $\mathcal Z_v$ independent of $A$, and preserve the
selection bounds of Lemma~\ref{lem:core} up to universal constants.

\begin{lemma}[Discretization of the retained weights]\label{lem:covering}
There are universal constants $C_0,C>0$ and finite families
$\mathcal Z_v$, $k\le v\le d$, depending only on the dimensions,
with the following properties. Each $z=(z_1,\ldots,z_k)\in\mathcal Z_v$
consists of nonnegative unit factors with supports
$I_j\subseteq[n_j]$ satisfying $\sum_j|I_j|=v$.
For $e\in\prod_jI_j$, write $z_e=\prod_jz_{j,e_j}$ and
$w_z(e)=\sum_jz_{j,e_j}^{-2}$. Every nonzero coordinate satisfies
\begin{equation}\label{eq:coordinate-floor}
 z_{j,i}^2\ge\frac1{C_0v\log^2(ed/v)}
 \quad(z\in\mathcal Z_v,\ i\in I_j,\ 1\le j\le k).
\end{equation}
The cardinality of each family satisfies
\begin{equation}\label{eq:entropy}
 |\mathcal Z_v|\le\exp\left\{Cv\log\left(\frac{ed}{v}\right)\log k\right\}.
\end{equation}

For $u,F,t>0$ and $A\in\{H_t(A)>u,\ f(A)\le F\}$, let
$x,E,S,\mu,J$ be as constructed in
Section~\ref{sec:entry-coordinate-selection}.
With $I_j=\pi_jJ$ and $v=\sum_j|I_j|$, there is a tuple
$z\in\mathcal Z_v$ supported on these $I_j$ for which the
selection bounds take the form
\begin{equation}\label{eq:rounded-selection}
 p_{j,i}(J)\ge\frac{c_1}{3}z_{j,i}^2\quad(i\in I_j),
 \qquad \mu_ew_z(e)\le8\mathcal E\quad(e\in J).
\end{equation}
The set $J$ still satisfies $\mu(J)>1/2$. Its weighted
contributions computed with $z$ satisfy
\begin{align}
 A_ez_e&>\frac t2,\qquad
 A_ez_ew_z(e)\le\frac{C_0kF^2}{S}\quad(e\in J),
 \label{eq:rounded-terms}\\
 \sum_{e\in J}A_ez_ew_z(e)&\ge\frac{c_1Sv}{6}.
 \label{eq:rounded-sum}
\end{align}
\end{lemma}
\begin{proof}
We first round and normalize the retained weights $p_{j,i}$. Next, we
transfer the selection bounds to $z$ and control the resulting
weighted contributions. Finally, we count the possible tuples $z$.
Fix $u,F,t>0$ and $A\in\{H_t(A)>u,\ f(A)\le F\}$.
Let $x,E,S,J$ and $I_j=\pi_jJ$ be as in
Section~\ref{sec:entry-coordinate-selection}, and put
$v=\sum_j|I_j|$. In particular, $S>u$ and
\eqref{eq:core-sum}--\eqref{eq:core-bound} hold.

\emph{Step 1. Rounding and normalization.}
Choose a common grid spacing $\delta$ and round each squared
coordinate $p_{j,i}$ upward by setting, for a universal
$C_0\ge16/c_2$,
\[
 \delta=\frac{2}{C_0kv\log^2(ed/v)},\qquad
 m_{j,i}=\left\lceil\frac{p_{j,i}}{\delta}\right\rceil
 \quad(i\in I_j).
\]
The rounded weight is $\delta m_{j,i}$.
Upward rounding and \eqref{eq:floor} give
\[
 0\le\frac{\delta m_{j,i}-p_{j,i}}{p_{j,i}}
 <\frac{\delta}{p_{j,i}}
 \le\frac{2}{C_0c_2k}\le\frac1{8k}.
\]
To bound the total rounded weight in each mode, use
$p_{j,i}\le\delta m_{j,i}\le p_{j,i}+\delta$ and
$p_j(I_j)\ge\mu(J)>1/2$ to obtain
\[
 \frac12<p_j(I_j)\le\delta\sum_{i\in I_j}m_{j,i}
 \le p_j(I_j)+|I_j|\delta\le1+v\delta.
\]
Here $v\delta=2/[C_0k\log^2(ed/v)]\le1/(8k)$.

We normalize the rounded weights so that $\sum_{i\in I_j}z_{j,i}^2=1$, the
identity used in the moment estimates of
Section~\ref{sec:probability}. Set
\begin{equation}\label{eq:rounding}
 z_{j,i}
 =\left(\frac{\delta m_{j,i}}{\sum_{a\in I_j}\delta m_{j,a}}\right)^{1/2}
 =\left(\frac{m_{j,i}}{\sum_{a\in I_j}m_{j,a}}\right)^{1/2}
 \quad(i\in I_j),
\end{equation}
and extend $z_j$ by zero outside $I_j$.
Dividing the bounds on each rounded weight by the bounds on its
mode total gives, for $i\in I_j$,
\begin{equation}\label{eq:square-comparison}
 \frac{p_{j,i}}{1+v\delta}\le z_{j,i}^2
 \le2\left(1+\frac1{8k}\right)p_{j,i}\le3p_{j,i},
 \qquad z_{j,i}^2\ge\frac1{C_0v\log^2(ed/v)}.
\end{equation}
The final lower bound uses \eqref{eq:floor} and $C_0c_2\ge16$.

\emph{Step 2. Transferring the selection and contribution bounds.}
We first prove \eqref{eq:rounded-selection}. The coordinate
bounds in \eqref{eq:square-comparison} imply
\[
 \frac1{3p_{j,i}}\le\frac1{z_{j,i}^2}
 \le\frac{1+v\delta}{p_{j,i}}\le\frac2{p_{j,i}}
 \quad(i\in I_j).
\]
Setting $i=e_j$ and summing over $j$ gives
\[
 \frac{w_e}{3}\le\sum_{j=1}^k\frac1{z_{j,e_j}^2}
 =w_z(e)\le2w_e\quad(e\in J).
\]
Together with $z_{j,i}^2\le3p_{j,i}$, \eqref{eq:retention}
and \eqref{eq:score}, this yields
\[
 p_{j,i}(J)\ge c_1p_{j,i}\ge\frac{c_1}{3}z_{j,i}^2,
 \qquad \mu_ew_z(e)\le2\mu_ew_e\le8\mathcal E.
\]
This proves \eqref{eq:rounded-selection}. The mass bound
$\mu(J)>1/2$ and the coordinate bound \eqref{eq:coordinate-floor}
complete the transfer of the four conclusions of Lemma~\ref{lem:core}.

To transfer the bounds on $A_ex_ew_e$ to $A_ez_ew_z(e)$, we
still need to compare $z_e$ with $x_e$, with constants independent
of $k$. By \eqref{eq:rounding}, for $e\in J$,
\[
 z_e=\left(
 \frac{\prod_{j=1}^k\delta m_{j,e_j}}
 {\prod_{j=1}^k\sum_{i\in I_j}\delta m_{j,i}}
 \right)^{1/2}.
\]
For the upper bound, Step~1 gives
$\delta m_{j,e_j}\le(1+1/(8k))p_{j,e_j}$ and
$\sum_{i\in I_j}\delta m_{j,i}\ge p_j(I_j)$.
Applying \eqref{eq:product} with $b_j=\one_{I_j}$ also gives
\[
 \frac12<\mu(J)\le\E_\mu\prod_j\one_{I_j}(X_j)
 \le\prod_jp_j(I_j)^{1/2}.
\]
Since $\prod_jp_{j,e_j}^{1/2}=x_e$, these inequalities imply
\[
 z_e\le\frac{(1+1/(8k))^{k/2}}{\prod_jp_j(I_j)^{1/2}}x_e
 \le\frac{(1+1/(8k))^{k/2}}{\mu(J)}x_e<3x_e,
\]
where the last inequality uses $\mu(J)>1/2$ and
$(1+1/(8k))^{k/2}\le\exp(1/16)<3/2$.
For the lower bound, use $\delta m_{j,e_j}\ge p_{j,e_j}$ and
$\sum_{i\in I_j}\delta m_{j,i}\le1+v\delta$ to obtain
\[
 z_e\ge(1+v\delta)^{-k/2}x_e\ge\frac{x_e}{2},
\]
where the last inequality follows from $v\delta\le1/(8k)$.
The lower bound $z_e\ge x_e/2$ and $A_ex_e>t$ give
$A_ez_e>t/2$. We now combine the comparisons of $z_e$ with
$x_e$ and of $w_z(e)$ with $w_e$ to obtain
\[
 \frac16A_ex_ew_e\le A_ez_ew_z(e)\le6A_ex_ew_e
 \quad(e\in J).
\]
The upper bound and \eqref{eq:core-bound} give
$A_ez_ew_z(e)\le24kF^2/S\le C_0kF^2/S$ for $C_0\ge24$,
proving \eqref{eq:rounded-terms}. Summing the lower bound
and applying \eqref{eq:core-sum} proves \eqref{eq:rounded-sum}.

\emph{Step 3. Constructing and counting the finite families.}
In \eqref{eq:rounding}, fixing the supports $I_j$ and the positive
integers $m_{j,i}$ fixes every coordinate of $z_j$. To obtain
finitely many tuples independent of $A$, we bound
$\sum_j\sum_{i\in I_j}m_{j,i}$. Since $p_j(I_j)\le1$ in each mode,
\begin{equation}\label{eq:integer-budget}
 \sum_j\sum_{i\in I_j}m_{j,i}
 \le\frac{k}{\delta}+v
 =\frac{C_0}{2}k^2v\log^2\!\left(\frac{ed}{v}\right)+v
 \le C_0k^2v\log^2\!\left(\frac{ed}{v}\right).
\end{equation}
For each $k\le v\le d$, let $I=(I_1,\ldots,I_k)$ range over
nonempty supports $I_j\subseteq[n_j]$, and let
$m=(m_{j,i})_{j,i}$ range over positive integer arrays satisfying
\begin{equation}\label{eq:family-parameters}
 \begin{gathered}
  \sum_{j=1}^k|I_j|=v,\qquad
  \sum_{j=1}^k\sum_{i\in I_j}m_{j,i}
  \le C_0k^2v\log^2\!\left(\frac{ed}{v}\right),\\
  \frac{m_{j,i}}{\sum_{a\in I_j}m_{j,a}}
  \ge\frac1{C_0v\log^2\!\left(\frac{ed}{v}\right)}
  \quad(i\in I_j,\ 1\le j\le k).
 \end{gathered}
\end{equation}
For each such pair $(I,m)$, define $z(I,m)=(z_1,\ldots,z_k)$ by
\[
 z_{j,i}=\begin{cases}
  \displaystyle\left(\frac{m_{j,i}}{\sum_{a\in I_j}m_{j,a}}\right)^{1/2},
    &i\in I_j,\\
  0,&i\notin I_j,
 \end{cases}
 \qquad 1\le j\le k,
\]
and set
\[
 \mathcal Z_v=\bigl\{z(I,m)\mid
 I,m\text{ satisfy }\eqref{eq:family-parameters}\bigr\}.
\]
The supports $I_j=\pi_jJ$ and integers
$m_{j,i}=\lceil p_{j,i}/\delta\rceil$ used in Step~1 satisfy
\eqref{eq:family-parameters} by \eqref{eq:square-comparison} and
\eqref{eq:integer-budget}. Hence the tuple $z$ defined in
\eqref{eq:rounding} belongs to $\mathcal Z_v$.

The coordinate bound \eqref{eq:coordinate-floor} follows from
\eqref{eq:family-parameters}. To prove \eqref{eq:entropy}, we count
the choices of $I$ and $m$. The supports
$I_j=\operatorname{supp}(z_j)$ together specify the $v$-element set
$\{(j,i):1\le j\le k,\ i\in I_j\}$.
Since there are $d=\sum_jn_j$ possible pairs $(j,i)$ with
$1\le j\le k$ and $i\in[n_j]$, there are at most $\binom dv$
choices of $I$. For each fixed $I$, order its $v$ coordinate
pairs and regard $m$ as a list of $v$ positive integers.
If their sum is $s\ge v$, choosing $v-1$ separators among the
$s-1$ gaps between $s$ units gives $\binom{s-1}{v-1}$ such lists.
Thus, for an integer $B\ge v$, the number of arrays $m$ with
$\sum_j\sum_{i\in I_j}m_{j,i}\le B$ is
\[
 \sum_{s=v}^B\binom{s-1}{v-1}
 =\sum_{s=v}^B\left[\binom sv-\binom{s-1}{v}\right]
 =\binom Bv,
\]
where the first equality uses Pascal's identity and the second
sum telescopes.
Taking $B=\lfloor C_0k^2v\log^2(ed/v)\rfloor$ gives
\[
 |\mathcal Z_v|\le\binom dv\binom Bv
 \le\binom dv\left(C_0e k^2\log^2\left(\frac{ed}{v}\right)\right)^v.
\]
Since $\binom dv\le\exp\{v\log(ed/v)\}$ and
$\log\log(ed/v)\le\log(ed/v)$, this implies
\[
 |\mathcal Z_v|
 \le\exp\left\{Cv\left[\log\left(\frac{ed}{v}\right)+\log k\right]\right\}
 \le\exp\left\{Cv\log\left(\frac{ed}{v}\right)\log k\right\}.
\]
The last step uses $\log(ed/v)\ge1$ and $\log k\ge\log2$,
after increasing $C$.
Thus $\mathcal Z_v$ is finite and satisfies \eqref{eq:entropy}.
\end{proof}

To obtain a finite covering of $\{H_t(A)>u,\ f(A)\le F\}$,
define, for each $z\in\mathcal Z_v$,
\begin{equation}\label{eq:statistic}
 T_z(A)=\sum_{e\in\prod_jI_j}
 \left(\min\left\{A_ez_ew_z(e),\frac{C_0kF^2}{u}\right\}
       -\frac t8w_z(e)\right)_+,
\end{equation}
where the minimum bounds every summand by $C_0kF^2/u$ over the
full product $\prod_jI_j$, since \eqref{eq:rounded-terms}
gives a bound only on $J$.

\begin{lemma}[Finite covering]\label{lem:finite-covering}
With $C_0$ chosen sufficiently large, for all $u,F,t>0$,
\begin{equation}\label{eq:inclusion}
 \{H_t(A)>u,\ f(A)\le F\}
 \subseteq\bigcup_{v=k}^d\ \bigcup_{z\in\mathcal Z_v}
 \{A\ge0:T_z(A)>uv/C_0\}.
\end{equation}
\end{lemma}
\begin{proof}
For $A\in\{H_t(A)>u,\ f(A)\le F\}$, take $v,z,J$ supplied
by Lemma~\ref{lem:covering}. For $e\in J$,
\eqref{eq:rounded-terms} and $S>u$ give
$A_ez_ew_z(e)<C_0kF^2/u$ and $A_ez_e>t/2$.
Restricting the nonnegative sum in \eqref{eq:statistic} to $J$
and applying \eqref{eq:rounded-sum} therefore gives
\[
 T_z(A)\ge\sum_{e\in J}\left(A_ez_e-\frac t8\right)w_z(e)
 \ge\frac34\sum_{e\in J}A_ez_ew_z(e)
 \ge\frac{c_1Sv}{8}>\frac{uv}{C_0},
\]
where $C_0\ge8/c_1$. Since $A$ was arbitrary, this proves
\eqref{eq:inclusion}.
\end{proof}

\subsection{Probability estimates}\label{sec:probability}
For a nonnegative random tensor $A$, Lemma~\ref{lem:finite-covering}
and the union bound give
\[
 \PP\{H_t(A)>u,\ f(A)\le F\}
 \le\sum_{v=k}^d\sum_{z\in\mathcal Z_v}
       \PP\{T_z(A)>uv/C_0\}.
\]
The cardinality bound \eqref{eq:entropy} controls the number of
terms in this union bound. Lemma~\ref{lem:statistic-tail}
estimates each probability for a fixed $z\in\mathcal Z_v$.
Throughout this subsection, $C_0$ is fixed as in
Lemma~\ref{lem:finite-covering}.

\begin{lemma}[Tail bound for $T_z$]\label{lem:statistic-tail}
There are universal constants $c,C>0$ with the following property.
Let $A=(A_e)_{e\in\mathcal I}$ have independent nonnegative entries
satisfying $\E A_e^2\le1/n$. Define $T_z$ by
\eqref{eq:statistic} with $t=1/(n\sqrt{k\log k})$.
For every $F\ge1$, $u\ge CF\sqrt{k\log k}$, $k\le v\le d$,
and deterministic $z\in\mathcal Z_v$,
\begin{equation}\label{eq:pointwise}
 \PP(T_z(A)>uv/C_0)
 \le\exp\left\{-c\frac{u^2}{kF^2}v\log(ed/v)\right\}.
\end{equation}
\end{lemma}
\begin{proof}
The constants $C,c>0$ below are universal and may change from
line to line. Fix $v$ and a deterministic $z\in\mathcal Z_v$.
For $e\in\prod_jI_j$, set
\[
 X_e=\left(\min\left\{A_ez_ew_z(e),\frac{C_0kF^2}{u}\right\}
                  -\frac{t}{8}w_z(e)\right)_+.
\]
Then $T_z(A)=\sum_eX_e$, and the $X_e$ are independent with
$0\le X_e\le b:=C_0kF^2/u$.
To apply Bennett's inequality, we first bound $\E T_z(A)$ and
$\sum_e\E X_e^2$.

For the mean, the definition of $X_e$ and $\E A_e^2\le1/n$ give
\[
 \E X_e\le w_z(e)\E\left(A_ez_e-\frac{t}{8}\right)_+
 \le\frac2t\E A_e^2\,z_e^2w_z(e)
 \le\frac2{tn}z_e^2w_z(e).
\]
Since $z_e^2=\prod_jz_{j,e_j}^2$ and
$\sum_{i\in I_j}z_{j,i}^2=1$, summing over $e\in\prod_jI_j$ yields
\[
 \E T_z(A)\le\frac2{tn}\sum_ez_e^2w_z(e)
 =\frac2{tn}\sum_{j=1}^k|I_j|
 =\frac{2v}{tn}=2v\sqrt{k\log k}.
\]
Taking the universal lower bound on $u$ large enough ensures
$uv/C_0-\E T_z(A)\ge uv/(2C_0)$.

For the second moments, $X_e\le A_ez_ew_z(e)$ gives
$\E X_e^2\le n^{-1}z_e^2w_z(e)^2$.
Summing over $e\in\prod_jI_j$ and expanding $w_z(e)^2$, we obtain
\begin{equation}\label{eq:second}
 \sum_e\E X_e^2
 \le\frac1n\sum_ez_e^2w_z(e)^2
 =\frac1n\left(\sum_{j=1}^k\sum_{i\in I_j}\frac1{z_{j,i}^2}
                  +2\sum_{j<\ell}|I_j||I_\ell|\right)
 \le\frac{Cv^2\log^2(ed/v)}{n}.
\end{equation}
Here \eqref{eq:coordinate-floor} bounds the first sum by
$C_0v^2\log^2(ed/v)$, and
$2\sum_{j<\ell}|I_j||I_\ell|\le v^2$.

Bennett's inequality \cite[Theorem~2.9]{BLM} and
\eqref{eq:second} give, for $s>0$,
\begin{equation}\label{eq:bennett}
 \PP(T_z(A)-\E T_z(A)>s)\le
 \exp\left\{-\frac{Cv^2\log^2(ed/v)}{nb^2}
 h\left(\frac{nbs}{Cv^2\log^2(ed/v)}\right)\right\},
\end{equation}
where $h(a)=(1+a)\log(1+a)-a$ for $a\ge0$.
By the inequality $h(a)\ge(a/2)\log(1+a)$ for $a\ge0$ and
the mean bound, \eqref{eq:bennett} with $s=uv/(2C_0)$ gives
\begin{equation}\label{eq:rate-intermediate}
 \PP(T_z(A)>uv/C_0)
 \le\exp\left\{-c\frac{u^2v}{kF^2}
       \log\left(1+\frac{ckF^2n}{v\log^2(ed/v)}\right)\right\}.
\end{equation}

Since $kn\ge d$, $F\ge1$, and
$\log^2(ed/v)\le C\sqrt{d/v}$ for a universal $C\ge c$,
\[
 \log\left(1+\frac{ckF^2n}{v\log^2(ed/v)}\right)
 \ge\log\left(1+\frac cC\sqrt{d/v}\right)
 \ge\frac cC\log\left(1+\sqrt{d/v}\right)
 \ge\frac{c}{4C}\log(ed/v),
\]
where the last two inequalities use concavity of $\log(1+x)$ and
$d/v\ge1$.
Combining this estimate with \eqref{eq:rate-intermediate} proves
\eqref{eq:pointwise}.
\end{proof}

\begin{proof}[Proof of Theorem~\ref{thm:prob}]
We first prove the tail estimate \eqref{eq:tail}.
For $F\ge1$ and $u\ge CF\sqrt{k\log k}$, with $C$ sufficiently
large, apply the union bound to the cover in
Lemma~\ref{lem:finite-covering}. Lemma~\ref{lem:statistic-tail}
and the cardinality bound \eqref{eq:entropy} then give
\begin{align*}
 \PP\{H_t(A)>u,\ f(A)\le F\}
 &\le\sum_{v=k}^d
 \exp\left\{-\left(c\frac{u^2}{kF^2}-C\log k\right)
                         v\log(ed/v)\right\}\\
 &\le\sum_{v=k}^d
 \exp\left\{-\frac{cu^2}{2kF^2}v\log(ed/v)\right\}\\
 &\le d\exp\left\{-\frac{cu^2}{2kF^2}\log d\right\}
 \le\exp\left\{-\frac{cu^2}{4kF^2}\log d\right\}.
\end{align*}
The last line uses that $v\log(ed/v)$ is increasing on $[1,d]$
and $cu^2/(2kF^2)\ge2$.
This proves \eqref{eq:tail} after adjusting the universal constant $c$.

We next prove the expectation bound \eqref{eq:soft-expectation}.
Choose a universal $C_1\ge C$ such that $cC_1^2(\log2)^2\ge1$, where
$c,C$ are the constants in \eqref{eq:tail}.
For $s\ge1$, set
$u=C_1s\sqrt{k\log k}$ and $F=\sqrt s$.
Then $u\ge CF\sqrt{k\log k}$. Splitting the event according to
whether $f(A)^2>s$ and applying \eqref{eq:tail} gives
\begin{align*}
 &\PP(H_t(A)>C_1s\sqrt{k\log k})\\
 &\quad\le\PP(f(A)^2>s)
   +\PP\{H_t(A)>C_1s\sqrt{k\log k},\ f(A)^2\le s\}\\
 &\quad\le\PP(f(A)^2>s)
   +\exp\{-cC_1^2s\log k\log d\}\\
 &\quad\le\PP(f(A)^2>s)+e^{-s}.
\end{align*}
Integrating this bound and bounding the probability by one on
$0\le s\le1$, we obtain
\begin{align*}
 \frac{\E H_t(A)}{C_1\sqrt{k\log k}}
 &=\int_0^\infty
     \PP(H_t(A)>C_1s\sqrt{k\log k})\,ds\\
 &\le1+\int_1^\infty\PP(f(A)^2>s)\,ds
       +\int_1^\infty e^{-s}\,ds
 \le1+\E[f(A)^2]+e^{-1}.
\end{align*}
Taking $C=C_1(1+e^{-1})$ proves \eqref{eq:soft-expectation}.
This choice also satisfies the threshold requirement in \eqref{eq:tail}.
\end{proof}

\section{Sharpness for Rademacher tensors}\label{sec:sharpness}
This section shows that the upper bound in
Theorem~\ref{thm:main} is optimal up to a universal constant by
considering Rademacher tensors on $[2]^k$. More precisely, we prove
Theorem~\ref{thm:sign}, which gives
$\E\norm{R^{(k)}}_{\inj}\asymp\sqrt{k\log k}$ while
$f(R^{(k)})=\sqrt2$. The lower bound is obtained by constructing
a family of product vectors with a large Gaussian maximum and
comparing the Gaussian and Rademacher maxima over this family.

\subsection{The upper bound}
The approximation \eqref{eq:net-approximation} reduces the upper
bound to a finite maximum of Rademacher contractions, estimated
in Lemma~\ref{lem:rademacher-upper}.

\begin{lemma}[Rademacher upper bound]\label{lem:rademacher-upper}
Let $R^{(k)}$ have independent Rademacher entries on $[2]^k$.
For every $k\ge2$,
\[
 \E\norm{R^{(k)}}_{\inj}\le C\sqrt{k\log k},
 \qquad f(R^{(k)})=\sqrt2,
\]
where $C>0$ is universal.
\end{lemma}
\begin{proof}
Each fiber consists of two signs, so $f(R^{(k)})=\sqrt2$.
As in the proof of Lemma~\ref{lem:soft-norm}, choose a
$1/(2k)$-net $\mathcal N_j$ of the unit circle with
$|\mathcal N_j|\le(1+4k)^2$ and put
$\mathcal N=\prod_{j=1}^k\mathcal N_j$.
By \eqref{eq:net-approximation},
\[
 |\mathcal N|\le(1+4k)^{2k},\qquad
 \norm{R^{(k)}}_{\inj}
 \le2\max_{x\in\mathcal N}|R^{(k)}(x_1,\ldots,x_k)|.
\]

For each $x\in\mathcal N$, independence gives
\[
 \E e^{\lambda R^{(k)}(x_1,\ldots,x_k)}
 =\prod_e\cosh(\lambda x_e)
 \le\exp\left(\frac{\lambda^2}{2}\sum_e x_e^2\right)
 =e^{\lambda^2/2},\qquad \lambda\in\mathbb R,
\]
where $\sum_e x_e^2=\prod_j\norm{x_j}_2^2=1$.
The maximal inequality \cite[Theorem~2.5]{BLM}, applied to the
$2|\mathcal N|$ variables $\pm R^{(k)}(x_1,\ldots,x_k)$ with
$x\in\mathcal N$, yields
\begin{align*}
 \E\norm{R^{(k)}}_{\inj}
 &\le2\inf_{\lambda>0}
   \left\{\frac{\log(2|\mathcal N|)}{\lambda}+\frac{\lambda}{2}\right\}
 =2\sqrt{2\log(2|\mathcal N|)}\\
 &\le2\sqrt{2\log\{2(1+4k)^{2k}\}}
 \le C\sqrt{k\log k}.
\end{align*}
\end{proof}

\subsection{Product vectors for the Gaussian lower bound}
We seek a finite packing in the set of unit product vectors
\[
 \mathcal T=\left\{\bigotimes_{j=1}^k u_j:
          u_j\in\mathbb R^2,\ \norm{u_j}_2=1\right\}.
\]
For any tensor $A$ on $[2]^k$ and any $x^{(1)},\ldots,x^{(M)}\in\mathcal T$,
\[
 \norm A_{\inj}=\sup_{x\in\mathcal T}|\langle A,x\rangle|
 \ge\max_{m\le M}\langle A,x^{(m)}\rangle.
\]
By Sudakov minoration, an $\varepsilon_0$-packing with
$\log M\asymp k\log k$ gives a Gaussian lower bound of order
$\sqrt{k\log k}$, where $\varepsilon_0>0$ is a universal constant.
Lemma~\ref{lem:gaussian-product-packing} constructs such a packing.
Here we bound each tensor coordinate in absolute value by $(3/4)^k$
to control the comparison error in
Lemma~\ref{lem:gaussian-rademacher-comparison}.

\begin{lemma}[Product vectors with small coordinates]\label{lem:gaussian-product-packing}
There are universal constants $c,C>0$ such that, for all sufficiently
large $k$, one can find $M\ge2$ product vectors
$x^{(m)}=\bigotimes_{j=1}^k x_j^{(m)}$, with $x_j^{(m)}\in\mathbb R^2$ and
$\norm{x_j^{(m)}}_2=1$, satisfying
\begin{equation}\label{eq:packing-size}
 ck\log k\le\log M\le Ck\log k,
 \qquad\max_{\substack{m\le M\\e\in[2]^k}}|x_e^{(m)}|\le\left(\frac34\right)^k,
\end{equation}
where $x_e^{(m)}=\prod_{j=1}^k x_{j,e_j}^{(m)}$ for $e\in[2]^k$.
If $g=(g_e)_{e\in[2]^k}$ has independent standard Gaussian entries, then
\begin{equation}\label{eq:gaussian-packing}
 \E\max_{m\le M}\langle g,x^{(m)}\rangle\ge c\sqrt{k\log k}.
\end{equation}
\end{lemma}
\begin{proof}
We first parametrize a subset of $\mathcal T$ with small coordinates.
We then construct a packing of its angle parameters and transfer it
to the product vectors. Sudakov minoration gives
\eqref{eq:gaussian-packing}.

\emph{Step 1. Unit factors with small coordinates.}
To keep both coordinates of each unit factor close to $1/\sqrt2$,
fix $\delta=1/20$ and restrict the angle parameters to
\[
 \mathcal Q=\left[\frac{\pi}{4}-\delta,\frac{\pi}{4}+\delta\right]^k
 \subset\mathbb R^k.
\]
For $\theta=(\theta_1,\ldots,\theta_k)\in\mathcal Q$, put
\[
 x_j(\theta_j)=(\cos\theta_j,\sin\theta_j),\qquad
 x(\theta)=\bigotimes_{j=1}^k x_j(\theta_j).
\]
Each $\theta_j$ is the angle of the $j$th unit factor. Hence
$x(\mathcal Q)\subset\mathcal T$, and the restriction
$\theta\in\mathcal Q$ gives
\[
 \max_{e\in[2]^k}|x_e(\theta)|
 \le\cos^k\left(\frac{\pi}{4}-\delta\right)
 \le\left(\frac{1+\delta}{\sqrt2}\right)^k
 <\left(\frac34\right)^k.
\]

\emph{Step 2. A packing of the parameter set.}
Let $\Theta\subset\mathcal Q$ be a maximal Euclidean
$\delta$-packing. By maximality, the balls of radius $\delta$
centered at $\Theta$ cover $\mathcal Q$. Comparing volumes gives
\[
 |\Theta|\ge\frac{\operatorname{vol}(\mathcal Q)}
                 {\delta^k\operatorname{vol}(B_2^k)}
             =\frac{2^k}{\operatorname{vol}(B_2^k)},
\]
where $B_2^k$ is the Euclidean unit ball in $\mathbb R^k$.
The Gaussian integral gives
\[
 e^{-\frac{k}{2}}k^{\frac{k}{2}}\operatorname{vol}(B_2^k)
 \le\int_{\sqrt k B_2^k}\exp\left(-\frac{\norm{x}_2^2}{2}\right)\,dx
 \le(2\pi)^{\frac{k}{2}}.
\]
Consequently, for all sufficiently large $k$,
\[
 \log|\Theta|\ge\frac{k}{2}\log k-Ck\ge\frac{k}{4}\log k.
\]
Select $M=\lfloor k^{k/4}\rfloor$ distinct points
$\theta^{(1)},\ldots,\theta^{(M)}\in\Theta$ and put $x^{(m)}=x(\theta^{(m)})$.
Then $\log M\asymp k\log k$, and the coordinate bound from
Step 1 gives \eqref{eq:packing-size}.

\emph{Step 3. The induced packing of product vectors.}
We now show that $\lVert x^{(m)}-x^{(\ell)}\rVert_2\ge\varepsilon_0$ for
$m\ne\ell$, where $\varepsilon_0>0$ is a universal constant.
Fix $m\ne\ell$. Step 2 gives $\lVert\theta^{(m)}-\theta^{(\ell)}\rVert_2\ge\delta$,
and $\theta^{(m)},\theta^{(\ell)}\in\mathcal Q$ imply
$|\theta_j^{(m)}-\theta_j^{(\ell)}|\le2\delta=1/10$.
Thus the inequality $\log\cos u\le-u^2/2$ gives
\begin{equation}\label{eq:packing-metric}
 \langle x^{(m)},x^{(\ell)}\rangle
 =\prod_j\cos(\theta_j^{(m)}-\theta_j^{(\ell)})
 \le\exp\left(-\frac{\norm{\theta^{(m)}-\theta^{(\ell)}}_2^2}{2}\right)
 \le e^{-\delta^2/2}.
\end{equation}
Set $\varepsilon_0=\sqrt{2(1-e^{-\delta^2/2})}>0$.
Using \eqref{eq:packing-metric} and $\lVert x^{(m)}\rVert_2=\lVert x^{(\ell)}\rVert_2=1$,
we obtain
\[
 \E|\langle g,x^{(m)}-x^{(\ell)}\rangle|^2
 =\norm{x^{(m)}-x^{(\ell)}}_2^2
 =2\bigl(1-\langle x^{(m)},x^{(\ell)}\rangle\bigr)
 \ge\varepsilon_0^2.
\]
Thus $\{x^{(1)},\ldots,x^{(M)}\}\subset x(\mathcal Q)\subset\mathcal T$
is an $\varepsilon_0$-packing in the canonical metric of the
Gaussian process $x\mapsto\langle g,x\rangle$.
Sudakov minoration \cite[Theorem~7.4.1]{HDP} yields
\[
 \E\max_{m\le M}\langle g,x^{(m)}\rangle
 \ge c\varepsilon_0\sqrt{\log M}\ge c\sqrt{k\log k},
\]
where the last inequality uses \eqref{eq:packing-size} and the fact
that $\varepsilon_0$ is a positive universal constant.
This fixed packing scale therefore suffices for \eqref{eq:gaussian-packing}.
\end{proof}

\subsection{Comparison with Rademacher entries}
We now transfer the Gaussian lower bound in
Lemma~\ref{lem:gaussian-product-packing} to Rademacher entries.
Lemma~\ref{lem:gaussian-rademacher-comparison} compares the two
expected maxima over the same family $x^{(1)},\ldots,x^{(M)}$, with an
error controlled by the coordinate bound in \eqref{eq:packing-size}.

\begin{lemma}[Gaussian--Rademacher comparison]\label{lem:gaussian-rademacher-comparison}
Take $x^{(1)},\ldots,x^{(M)}$ from
Lemma~\ref{lem:gaussian-product-packing}.
Let $g=(g_e)_{e\in[2]^k}$ and
$\varepsilon=(\varepsilon_e)_{e\in[2]^k}$ be independent vectors
whose entries are independent standard Gaussian and Rademacher
variables, respectively. Then
\begin{equation}\label{eq:replacement}
 \left|\E\max_{m\le M}\langle\varepsilon,x^{(m)}\rangle
       -\E\max_{m\le M}\langle g,x^{(m)}\rangle\right|
 \le C(k\log k)^{3/4}\left(\frac{81}{128}\right)^{k/4},
\end{equation}
where $C>0$ is universal.
\end{lemma}
\begin{proof}
We use the smooth maximum from
\cite[proof of Theorem~1.4]{Chatterjee05} and a Lindeberg replacement
argument as in \cite[Section~2]{Chatterjee}, retaining the
fourth-order remainder. For $\beta>0$ and $z\in\mathbb R^{2^k}$, define
\[
 F_\beta(z)=\frac1\beta\log\sum_{m=1}^M
                         e^{\beta\langle z,x^{(m)}\rangle}.
\]
Then
\[
 0\le F_\beta(z)-\max_m\langle z,x^{(m)}\rangle\le\frac{\log M}{\beta}.
\]
To bound the fourth derivatives, put
\[
 \pi_m(z)=\frac{e^{\beta\langle z,x^{(m)}\rangle}}
              {\sum_{\ell=1}^M e^{\beta\langle z,x^{(\ell)}\rangle}},
 \qquad \bar x_e(z)=\sum_{m=1}^M\pi_m(z) x_e^{(m)}.
\]
Differentiation gives
\[
 \partial_e^4F_\beta(z)
 =\beta^3\left\{
   \sum_m\pi_m(z)(x_e^{(m)}-\bar x_e(z))^4
   -3\left[\sum_m\pi_m(z)(x_e^{(m)}-\bar x_e(z))^2\right]^2
             \right\}.
\]
Since $\sum_m\pi_m(z)=1$, the coordinate bound \eqref{eq:packing-size}
gives $|x_e^{(m)}-\bar x_e(z)|\le2(3/4)^k$, and hence
\begin{equation}\label{eq:fourth-derivative}
 |\partial_e^4F_\beta(z)|
 \le64\beta^3\left(\frac34\right)^{4k}
\end{equation}
uniformly in $z$.

Fix any ordering of $[2]^k$ and define
\[
 F_{\beta,e}(s)
 =F_\beta\bigl((\varepsilon_f)_{f<e},s,(g_f)_{f>e}\bigr),
 \qquad s\in\mathbb R.
\]
Taylor's formula and \eqref{eq:fourth-derivative} give
\[
 \left|F_{\beta,e}(s)
       -\sum_{r=0}^3\frac{s^r}{r!}F_{\beta,e}^{(r)}(0)\right|
 \le\frac{64\beta^3}{4!}\left(\frac34\right)^{4k}|s|^4.
\]
The random function $F_{\beta,e}$ is independent of
$\varepsilon_e$ and $g_e$. Their first three moments agree, so
\[
 \E\bigl[F_{\beta,e}^{(r)}(0)(\varepsilon_e^r-g_e^r)\bigr]
 =\E F_{\beta,e}^{(r)}(0)
       \bigl(\E\varepsilon_e^r-\E g_e^r\bigr)=0,
 \qquad r=1,2,3.
\]
Since $\E|\varepsilon_e|^4+\E|g_e|^4=1+3=4$, telescoping over
the entries yields
\begin{align*}
 |\E F_\beta(\varepsilon)-\E F_\beta(g)|
 &=\left|\sum_e\E\bigl[F_{\beta,e}(\varepsilon_e)
                       -F_{\beta,e}(g_e)\bigr]\right|\\
 &\le\frac{64\beta^3}{4!}\left(\frac34\right)^{4k}
       \sum_e\bigl(\E|\varepsilon_e|^4+\E|g_e|^4\bigr)
 =\frac{32}{3}\beta^3\left(\frac{81}{128}\right)^k.
\end{align*}
Choose $\beta=[(128/81)^k\log M]^{1/4}$ to balance the smoothing
and replacement errors. By \eqref{eq:packing-size},
\begin{align*}
 \left|\E\max_m\langle\varepsilon,x^{(m)}\rangle
       -\E\max_m\langle g,x^{(m)}\rangle\right|
 &\le\frac{2\log M}{\beta}
       +\frac{32}{3}\beta^3\left(\frac{81}{128}\right)^k\\
 &\le C(\log M)^{3/4}\left(\frac{81}{128}\right)^{k/4}
 \le C(k\log k)^{3/4}\left(\frac{81}{128}\right)^{k/4},
\end{align*}
which is \eqref{eq:replacement}.
\end{proof}

\begin{proof}[Proof of Theorem~\ref{thm:sign}]
Lemma~\ref{lem:rademacher-upper} gives the upper bound and
$f(R^{(k)})=\sqrt2$. For the lower bound, take $k$ sufficiently
large, choose the product vectors from
Lemma~\ref{lem:gaussian-product-packing}, and let
$\varepsilon$ be the entry vector of $R^{(k)}$.
Since each factor has unit norm,
\eqref{eq:gaussian-packing} and
Lemma~\ref{lem:gaussian-rademacher-comparison} give
\[
 \E\norm{R^{(k)}}_{\inj}
 \ge\E\max_{m\le M}\langle\varepsilon,x^{(m)}\rangle
 \ge c\sqrt{k\log k}
       -C(k\log k)^{3/4}\left(\frac{81}{128}\right)^{k/4}
 \ge\frac c2\sqrt{k\log k}.
\]
For the remaining finitely many orders, use
\[
 \norm{R^{(k)}}_{\inj}\ge f(R^{(k)})=\sqrt2
\]
and decrease the universal lower constant in \eqref{eq:sign} if necessary.
\end{proof}

\section{Applications and the role of the hypotheses}\label{sec:applications}
We first prove the higher-moment extension and explain the roles of
the largest entry and the moment order. We then apply the expectation
comparison to sparse Bernoulli tensors and give a matrix example
showing why identical distribution cannot be omitted in general.

\subsection{Higher moments}\label{sec:moments}
The expectation comparison in Theorem~\ref{thm:main} supplies the
dependence on the tensor order. To prove
Corollary~\ref{cor:higher-moments}, we control the additional
$L_q$ contribution by the largest entry. This moment estimate
requires only independence and centering.

\begin{proof}[Proof of Corollary~\ref{cor:higher-moments}]
Put $Z=\norm T_{\inj}$ and $M=M(T)$.
For independent centered entries with finite $q$th absolute moments,
we first prove
\begin{equation}\label{eq:norm-moment-increment}
 \norm Z_{L_q}\le C\bigl(\E Z+\sqrt q\,\norm M_{L_q}\bigr),
 \qquad q\ge1.
\end{equation}
We prove this estimate separately for $q\ge2$ and $1\le q<2$,
then apply Theorem~\ref{thm:main}.

\emph{Step 1. Moments of order $q\ge2$.}
Let $T^{(e)}$ replace $T_e$ by an independent copy $T_e'$, and put
\[
 Z^{(e)}=\norm{T^{(e)}}_{\inj},\qquad
 V_+=\sum_e\E\bigl[(Z-Z^{(e)})_+^2\mid T\bigr].
\]
The moment Efron--Stein inequality \cite[Theorem~2]{BBLM} gives
\[
 \norm{(Z-\E Z)_+}_{L_q}
 \le C\sqrt q\,\norm{\sqrt{V_+}}_{L_q}.
\]
To bound $V_+$, choose $a=a_1\otimes\cdots\otimes a_k$ with
$\norm{a_j}_2=1$ and $Z=\sum_e a_eT_e$.
Then $\sum_e a_e^2=1$. The supremum defining $Z^{(e)}$ includes
the choice $x_j=a_j$, so
\begin{align*}
 Z^{(e)}
 &=\sup_{\norm{x_j}_2\le1}
   \left|\left\langle T^{(e)},x_1\otimes\cdots\otimes x_k\right\rangle\right|
 \ge\left|\left\langle T^{(e)},a\right\rangle\right|\\
 &\ge\sum_{f\ne e}a_fT_f+a_eT_e'
 =Z+a_e(T_e'-T_e).
\end{align*}
Rearranging and taking positive parts gives
\[
 (Z-Z^{(e)})_+
 \le\bigl(a_e(T_e-T_e')\bigr)_+
 \le |a_e|\,|T_e-T_e'|.
\]
Independence and centering therefore give
\[
 V_+\le\sum_e a_e^2\E\bigl[(T_e-T_e')^2\mid T\bigr]
     =\sum_e a_e^2(T_e^2+\E T_e^2)
     \le M^2+\E M^2.
\]
Since $q\ge2$,
\[
 \norm{\sqrt{V_+}}_{L_q}
 \le\norm M_{L_q}+\norm M_{L_2}
 \le2\norm M_{L_q}.
\]
Since $Z\le\E Z+(Z-\E Z)_+$, Minkowski's inequality gives
\[
 \norm Z_{L_q}
 \le\E Z+\norm{(Z-\E Z)_+}_{L_q}
 \le\E Z+C\sqrt q\,\norm M_{L_q},
\]
which proves \eqref{eq:norm-moment-increment} for $q\ge2$.

\emph{Step 2. Moments of order $1\le q<2$.}
Let $T'$ be an independent copy of $T$. Since
$\E[T-T'\mid T]=T$, conditional Jensen's inequality gives
\[
 \norm Z_{L_q}\le\bigl\|\norm{T-T'}_{\inj}\bigr\|_{L_q}.
\]
Put $N=|\mathcal I|$ and enumerate the tensor coordinates by
$i=1,\ldots,N$.
Write $T_i,T_i'$ for the corresponding entries of $T,T'$, and
let $E_i$ have a one at coordinate $i$ and zeros elsewhere. Then
\[
 T-T'=\sum_{i=1}^N V_i,
\]
where $V_i=(T_i-T_i')E_i$ for $1\le i\le N$.
These tensors are independent and symmetric.
To apply the moment inequality for maximal partial sums in
\cite[Theorem~6.2]{HMS}, define
\[
 U=\max_{1\le j\le N}\Bigl\|\sum_{i=1}^jV_i\Bigr\|_{\inj},
 \qquad
 W=\max_{1\le i\le N}\norm{V_i}_{\inj}
  =\max_{1\le i\le N}|T_i-T_i'|.
\]
This inequality gives
\[
 \norm U_{L_q}\le C\bigl(\E U+\norm W_{L_q}\bigr),
 \qquad 1\le q<2.
\]
Since the $V_i$ are independent and symmetric, L\'evy's maximal
inequality \cite[equation~(5)]{BT} gives
\[
 \PP(U>t)\le2\PP(\norm{T-T'}_{\inj}>t),\qquad t>0.
\]
Integrating over $t>0$ and using the triangle inequality, we obtain
\begin{align*}
 \E U
 &=\int_0^\infty\PP(U>t)\,dt
 \le2\int_0^\infty\PP(\norm{T-T'}_{\inj}>t)\,dt\\
 &=2\E\norm{T-T'}_{\inj}
 \le2\bigl(\E\norm T_{\inj}+\E\norm{T'}_{\inj}\bigr)
 =4\E Z.
\end{align*}
To bound $\norm W_{L_q}$, we use
\[
 W\le M(T)+M(T'),\qquad \norm W_{L_q}\le2\norm M_{L_q}.
\]
Since $\norm{T-T'}_{\inj}\le U$, we obtain
\[
 \norm Z_{L_q}
 \le\norm U_{L_q}
 \le C\bigl(\E Z+\norm M_{L_q}\bigr),
\]
which proves \eqref{eq:norm-moment-increment} for $1\le q<2$.

\emph{Step 3. The fiber comparison.}
For i.i.d.\ entries, Theorem~\ref{thm:main} and
\eqref{eq:norm-moment-increment} give
\[
 \norm Z_{L_q}
 \le C\left(\sqrt{k\log k}\,\E f(T)
              +\sqrt q\,\norm{M(T)}_{L_q}\right),
\]
which is \eqref{eq:moment-refined}.
Since $M(T)\le f(T)\le Z$ and $\E f(T)\le\norm{f(T)}_{L_q}$,
\eqref{eq:moment-comparison} follows.
\end{proof}

\begin{remark}[The largest-entry term]
The largest-entry term in \eqref{eq:moment-refined} gives a sharper
bound for Rademacher tensors. For such a tensor $R$ on $[n]^k$,
\[
 M(R)=1,\qquad f(R)=\sqrt n,\qquad
 \bigl\|\norm R_{\inj}\bigr\|_{L_q}
 \le C\bigl(\sqrt{nk\log k}+\sqrt q\bigr).
\]
Here the term depending on $q$ is $\sqrt q$, while
\eqref{eq:moment-comparison} gives $\sqrt{nq}$.
Keeping $M(R)$ therefore saves a factor $\sqrt n$ in this term.
\end{remark}

\begin{remark}[Sharpness of the moment term]
The factor $\sqrt q$ cannot be improved uniformly over the tensor
order. To see this, take $R^{(k)}$ from Theorem~\ref{thm:sign} and
put $q_k=2^k$.
On the event that all $q_k$ entries equal one, whose probability is
$2^{-q_k}$, the injective norm is $\sqrt{q_k}$. Hence
\[
 \frac{\bigl\|\lVert R^{(k)}\rVert_{\inj}\bigr\|_{L_{q_k}}}
      {\bigl\|f(R^{(k)})\bigr\|_{L_{q_k}}}
 \ge\frac{\sqrt{q_k}(2^{-q_k})^{1/q_k}}{\sqrt2}
 =\frac{\sqrt{q_k}}{2\sqrt2}.
\]
Because $\sqrt{k\log k}=o(\sqrt{q_k})$, no bound of the form
\[
 \bigl\|\norm T_{\inj}\bigr\|_{L_q}
 \le C\bigl(\sqrt{k\log k}+h(q)\bigr)\norm{f(T)}_{L_q}
\]
can hold universally with $h(q)=o(\sqrt q)$.
This argument uses growing tensor order. For $n_1\times n_2$ matrices
($k=2$), Seginer's comparison has a universal constant for
$1\le q\le2\log\max\{n_1,n_2\}$ \cite[Corollary 2.2]{Seginer},
a range of uniformity not recovered by \eqref{eq:moment-comparison}.
\end{remark}

\subsection{Centered Bernoulli entries with parameter \texorpdfstring{$1/n$}{1/n}}
At parameter $1/n$, each fiber contains one nonzero entry on average
before centering. The squared norm of a centered fiber is determined
by its number of ones before centering. We estimate the largest
fiber count and then use
Theorem~\ref{thm:main} to obtain the expected injective norm.

\begin{corollary}\label{cor:bernoulli}
Fix $k\ge2$. Let $B$ have independent Bernoulli$(1/n)$ entries on
$[n]^k$, and put $T=B-\E B$. Then, as $n\to\infty$,
\begin{equation}\label{eq:bernoulli}
               \E\norm T_{\inj}\asymp_k\Phi(T)
               \asymp_k\sqrt{\frac{\log n}{\log\log n}}.
\end{equation}
\end{corollary}
\begin{proof}
We first express $f(T)$ in terms of the largest fiber count, then
bound that count from above using all fibers and from below using
an independent family of fibers.

\emph{Step 1. Reduction to fiber counts.}
By Theorem~\ref{thm:main}, it suffices to estimate $\Phi(T)=\E f(T)$.
Let $I\subset[n]^k$ be the index set of a fiber, and put
$Z_I=\sum_{e\in I}B_e$. Since $T_e=B_e-1/n$,
\[
 \sum_{e\in I}T_e^2
 =Z_I\left(1-\frac1n\right)^2+\frac{n-Z_I}{n^2}
 =\left(1-\frac2n\right)Z_I+\frac1n.
\]
Taking the maximum over all fibers and writing $D=\max_I Z_I$,
we obtain, for $n\ge4$,
\begin{equation}\label{eq:degree-fiber}
 f(T)^2=\left(1-\frac2n\right)D+\frac1n,\qquad
 \frac D2\le f(T)^2\le D+\frac1n.
\end{equation}

\emph{Step 2. Upper bound.}
Each fiber count has the distribution
$Z\sim\operatorname{Binomial}(n,1/n)$.
For every integer $m\ge1$, the variable $\binom Zm$ counts the
$m$-element subsets of successful trials and is at least one on
$\{Z\ge m\}$. Independence therefore gives
\[
 \PP(Z\ge m)\le\E\binom Zm=\binom nm n^{-m}
 \le\left(\frac em\right)^m.
\]
A union bound over the $kn^{k-1}$ fibers now yields
\begin{equation}\label{eq:bernoulli-degree-tail}
 \PP(D\ge m)\le kn^{k-1}\left(\frac em\right)^m.
\end{equation}
Choose
\[
 m_+=\left\lceil\frac{(k+2)\log n}{\log\log n}\right\rceil.
\]
Then, as $n\to\infty$,
\[
 \frac{m_+\log\left(\frac{m_+}{e}\right)}{\log n}
 \longrightarrow k+2.
\]
Thus, for all sufficiently large $n$,
\[
 kn^{k-1}\left(\frac e{m_+}\right)^{m_+}\le\frac1n.
\]
Also, $m_+\to\infty$, so $e/m_+<1/2$ for all sufficiently
large $n$. Using $\E D=\sum_{m\ge1}\PP(D\ge m)$ and
\eqref{eq:bernoulli-degree-tail}, we obtain
\[
 \E D\le m_++kn^{k-1}\sum_{m=m_+}^{\infty}
                \left(\frac e{m_+}\right)^m
      \le m_++2.
\]
Jensen's inequality and \eqref{eq:degree-fiber} then imply
\[
 \Phi(T)\le\sqrt{\E f(T)^2}
      \le\sqrt{m_++2+\frac1n}
      \le C_k\sqrt{\frac{\log n}{\log\log n}}.
\]

\emph{Step 3. Lower bound.}
For each of the $n^{k-1}$ choices of
$(i_2,\ldots,i_k)\in[n]^{k-1}$, define the corresponding fiber count by
\[
 Z_{i_2,\ldots,i_k}=\sum_{i_1=1}^n B_{i_1,i_2,\ldots,i_k}.
\]
The sums for distinct choices of $(i_2,\ldots,i_k)$ involve
disjoint sets of independent entries, so these counts are
independent copies of $Z$.
Since $\binom nm\ge(n/m)^m$,
\[
 \PP(Z\ge m)\ge\PP(Z=m)
 \ge m^{-m}\left(1-\frac1n\right)^n\ge\frac14m^{-m},
 \qquad 1\le m\le n,\quad n\ge2.
\]
Choose
\[
 m_-=\left\lfloor\frac{(k-1)\log n}{2\log\log n}\right\rfloor.
\]
For all sufficiently large $n$, we have $1\le m_-\le n$ and
$m_-\log m_-\le(3/4)(k-1)\log n$.
Since $D$ is the maximum over all fiber counts,
\[
 \{D<m_-\}\subseteq
 \bigcap_{(i_2,\ldots,i_k)\in[n]^{k-1}}
 \{Z_{i_2,\ldots,i_k}<m_-\}.
\]
Independence of these counts therefore gives
\begin{align*}
 \PP(D<m_-)
 &\le\bigl(1-\PP(Z\ge m_-)\bigr)^{n^{k-1}}
 \le\exp\left(-\frac14 n^{k-1}m_-^{-m_-}\right)\\
 &\le\exp\left(-\frac14 n^{\frac{k-1}{4}}\right)
 \longrightarrow0.
\end{align*}
On $\{D\ge m_-\}$, \eqref{eq:degree-fiber} gives
$f(T)\ge\sqrt{m_-/2}$. Using also $f(T)\ge0$, we obtain
\[
 \Phi(T)=\E f(T)
 \ge\E\!\left[f(T)\mathbf1_{\{D\ge m_-\}}\right]
 \ge\sqrt{\frac{m_-}{2}}\,\PP(D\ge m_-)
       \ge c_k\sqrt{\frac{\log n}{\log\log n}}
\]
for all sufficiently large $n$. Combining the two bounds on
$\Phi(T)$ with Theorem~\ref{thm:main} proves \eqref{eq:bernoulli}.
\end{proof}

\subsection{Why identical distribution matters}
Theorem~\ref{thm:main} assumes that the entries have a common
distribution. Proposition~\ref{prop:non-iid} shows that independence
and centering alone do not yield a comparison constant uniform
over the matrix dimensions.

\begin{proposition}\label{prop:non-iid}
There is a sequence of centered random matrices with independent
entries such that
\[
                       \frac{\E\norm A_{\inj}}{\E f(A)}
                       \longrightarrow\infty.
\]
\end{proposition}
\begin{proof}
Fix an integer $m\ge1$, and set
\[
 M=2^{m^2},\qquad A=\operatorname{diag}(R_1,\ldots,R_M),
\]
where $R_1,\ldots,R_M$ are independent $m\times m$ matrices with
i.i.d.\ Rademacher entries. Then $A$ has independent centered
entries, and its block structure gives
\[
 f(A)=\max_{1\le\ell\le M}f(R_\ell)=\sqrt m,\qquad
 \norm A_{\inj}=\max_{1\le\ell\le M}\norm{R_\ell}_{\inj}.
\]
Since the all-ones matrix $J_m=(1)_{i,j=1}^m$ has
$\norm{J_m}_{\inj}=m$,
\[
 \bigcup_{\ell=1}^M\{R_\ell=J_m\}
 \subseteq\{\norm A_{\inj}\ge m\}.
\]
For each $\ell$, independence of the $m^2$ entries gives
\[
 \PP(R_\ell=J_m)=\left(\frac12\right)^{m^2}=\frac1M.
\]
Taking probabilities in the inclusion and using independence of
the blocks yields
\[
 \PP\{\norm A_{\inj}\ge m\}
 \ge\PP\left(\bigcup_{\ell=1}^M\{R_\ell=J_m\}\right)
 =1-\left(1-\frac1M\right)^M\ge1-e^{-1}.
\]
Consequently,
\[
 \frac{\E\norm A_{\inj}}{\E f(A)}
 \ge\frac{m\,\PP\{\norm A_{\inj}\ge m\}}{\sqrt m}
 \ge(1-e^{-1})\sqrt m\longrightarrow\infty
 \qquad\text{as }m\to\infty.\qedhere
\]
\end{proof}

\end{document}